\documentclass[12pt,a4paper]{article}

\usepackage{amsmath,amssymb,latexsym,pstricks,mathrsfs,comment,amsthm,graphicx,tikz,enumerate,
pgffor,cite,wrapfig,multicol,float}
\usepackage{bibspacing}
\newcommand{\nc}{\newcommand}
\nc{\rnc}{\renewcommand}

\usepackage{a4wide}   \rnc{\ss}{} \nc{\ms}{} \nc{\bs}{}			

\makeatletter

\begin{document}

\hyphenation{mon-oid mon-oids}

\nc{\rasgp}{\ra_{\sgp}}
\nc{\itemit}[1]{\item[\emph{(#1)}]}
\nc{\E}{\mathcal E}
\nc{\TX}{\T(X)}
\nc{\TXP}{\T(X,\P)}
\nc{\EX}{\E(X)}
\nc{\EXP}{\E(X,\P)}
\nc{\SX}{\S(X)}
\nc{\SXP}{\S(X,\P)}
%\nc{\Sing}{\operatorname{Sing}}
\nc{\Sing}{\E}
\nc{\idrank}{\operatorname{idrank}}
\nc{\SingXP}{\Sing(X,\P)}
\nc{\De}{\Delta}
\nc{\sgp}{\operatorname{sgp}}
\nc{\mon}{\operatorname{mon}}
\nc{\Dn}{\mathcal D_n}
\nc{\Dm}{\mathcal D_m}

\nc{\lline}[1]{\draw(3*#1,0)--(3*#1+2,0);}
\nc{\uline}[1]{\draw(3*#1,5)--(3*#1+2,5);}
\nc{\thickline}[2]{\draw(3*#1,5)--(3*#2,0); \draw(3*#1+2,5)--(3*#2+2,0) ;}
\nc{\thicklabel}[3]{\draw(3*#1+1+3*#2*0.15-3*#1*0.15,4.25)node{{\tiny $#3$}};}

\nc{\slline}[3]{\draw(3*#1+#3,0+#2)--(3*#1+2+#3,0+#2);}
\nc{\suline}[3]{\draw(3*#1+#3,5+#2)--(3*#1+2+#3,5+#2);}
\nc{\sthickline}[4]{\draw(3*#1+#4,5+#3)--(3*#2+#4,0+#3); \draw(3*#1+2+#4,5+#3)--(3*#2+2+#4,0+#3) ;}
\nc{\sthicklabel}[5]{\draw(3*#1+1+3*#2*0.15-3*#1*0.15+#5,4.25+#4)node{{\tiny $#3$}};}

\nc{\stll}[5]{\sthickline{#1}{#2}{#4}{#5} \sthicklabel{#1}{#2}{#3}{#4}{#5}}
\nc{\tll}[3]{\stll{#1}{#2}{#3}00}

\nc{\mfourpic}[9]{
\slline1{#9}0
\slline3{#9}0
\slline4{#9}0
\slline5{#9}0
\suline1{#9}0
\suline3{#9}0
\suline4{#9}0
\suline5{#9}0
\stll1{#1}{#5}{#9}{0}
\stll3{#2}{#6}{#9}{0}
\stll4{#3}{#7}{#9}{0}
\stll5{#4}{#8}{#9}{0}
\draw[dotted](6,0+#9)--(8,0+#9);
\draw[dotted](6,5+#9)--(8,5+#9);
}
\nc{\vdotted}[1]{
\draw[dotted](3*#1,10)--(3*#1,15);
\draw[dotted](3*#1+2,10)--(3*#1+2,15);
}

\nc{\Clab}[2]{
\sthicklabel{#1}{#1}{{}_{\phantom{#1}}C_{#1}}{1.25+5*#2}0
}
\nc{\sClab}[3]{
\sthicklabel{#1}{#1}{{}_{\phantom{#1}}C_{#1}}{1.25+5*#2}{#3}
}
\nc{\Clabl}[3]{
\sthicklabel{#1}{#1}{{}_{\phantom{#3}}C_{#3}}{1.25+5*#2}0
}
\nc{\sClabl}[4]{
\sthicklabel{#1}{#1}{{}_{\phantom{#4}}C_{#4}}{1.25+5*#2}{#3}
}
\nc{\Clabll}[3]{
\sthicklabel{#1}{#1}{C_{#3}}{1.25+5*#2}0
}
\nc{\sClabll}[4]{
\sthicklabel{#1}{#1}{C_{#3}}{1.25+5*#2}{#3}
}

\nc{\mtwopic}[6]{
\slline1{#6*5}{#5}
\slline2{#6*5}{#5}
\suline1{#6*5}{#5}
\suline2{#6*5}{#5}
\stll1{#1}{#3}{#6*5}{#5}
\stll2{#2}{#4}{#6*5}{#5}
}
\nc{\mtwopicl}[6]{
\slline1{#6*5}{#5}
\slline2{#6*5}{#5}
\suline1{#6*5}{#5}
\suline2{#6*5}{#5}
\stll1{#1}{#3}{#6*5}{#5}
\stll2{#2}{#4}{#6*5}{#5}
\sClabl1{#6}{#5}{i}
\sClabl2{#6}{#5}{j}
}

%\nc{\keru}{\ker_{\operatorname{u}}} \nc{\kerl}{\ker_{\operatorname{l}}}
%\nc{\dom}{\operatorname{dom}_u} \nc{\codom}{\operatorname{dom}_l}
\nc{\keru}{\operatorname{ker}^\wedge} \nc{\kerl}{\operatorname{ker}_\vee}%\nc{\codom}{\operatorname{im}}

\nc{\coker}{\operatorname{coker}}
%\nc{\KER}{\operatorname{KER}}
\nc{\KER}{\ker}
\nc{\N}{\mathbb N}
\nc{\LaBn}{L_\al(\B_n)}
\nc{\RaBn}{R_\al(\B_n)}
\nc{\LaPBn}{L_\al(\PB_n)}
\nc{\RaPBn}{R_\al(\PB_n)}
\nc{\rhorBn}{\rho_r(\B_n)}
\nc{\DrBn}{D_r(\B_n)}
\nc{\DrPn}{D_r(\P_n)}
\nc{\DrPBn}{D_r(\PB_n)}
\nc{\DrKn}{D_r(\K_n)}
\nc{\alb}{\al_{\vee}}
\nc{\beb}{\be^{\wedge}}
\nc{\bnf}{\bn^\flat}
\nc{\Bal}{\operatorname{Bal}}
\nc{\Red}{\operatorname{Red}}
\nc{\Pnxi}{\P_n^\xi}
\nc{\Bnxi}{\B_n^\xi}
\nc{\PBnxi}{\PB_n^\xi}
\nc{\Knxi}{\K_n^\xi}
\nc{\C}{\mathscr C}
\nc{\exi}{e^\xi}
\nc{\Exi}{E^\xi}
\nc{\eximu}{e^\xi_\mu}
\nc{\Eximu}{E^\xi_\mu}
\nc{\REF}{ {\red [Ref?]} }
\nc{\GL}{\operatorname{GL}}
\rnc{\O}{\operatorname{O}}

\nc{\vtx}[2]{\fill (#1,#2)circle(.2);}
\nc{\lvtx}[2]{\fill (#1,0)circle(.2);}
\nc{\uvtx}[2]{\fill (#1,1.5)circle(.2);}

\nc{\Eq}{\mathfrak{Eq}}
%\nc{\Gau}{\Ga_{\operatorname{u}}} \nc{\Gal}{\Ga_{\operatorname{l}}}
\nc{\Gau}{\Ga^\wedge} \nc{\Gal}{\Ga_\vee}
%\nc{\Lamu}{\Lam_{\operatorname{u}}} \nc{\Laml}{\Lam_{\operatorname{l}}}
\nc{\Lamu}{\Lam^\wedge} \nc{\Laml}{\Lam_\vee}
\nc{\bX}{{\bf X}}
\nc{\bY}{{\bf Y}}
\nc{\ds}{\displaystyle}

\nc{\uvert}[1]{\fill (#1,1.5)circle(.2);}
\nc{\uuvert}[1]{\fill (#1,3)circle(.2);}
\nc{\uuuvert}[1]{\fill (#1,4.5)circle(.2);}
\rnc{\lvert}[1]{\fill (#1,0)circle(.2);}
\nc{\overt}[1]{\fill (#1,0)circle(.1);}
\nc{\overtl}[3]{\node[vertex] (#3) at (#1,0) {  {\tiny $#2$} };}
\nc{\cv}[2]{\draw(#1,1.5) to [out=270,in=90] (#2,0);}
\nc{\cvs}[2]{\draw(#1,1.5) to [out=270+30,in=90+30] (#2,0);}
\nc{\ucv}[2]{\draw(#1,3) to [out=270,in=90] (#2,1.5);}
\nc{\uucv}[2]{\draw(#1,4.5) to [out=270,in=90] (#2,3);}
\nc{\textpartn}[1]{{\lower0.45 ex\hbox{\begin{tikzpicture}[xscale=.2,yscale=0.2] #1 \end{tikzpicture}}}}
\nc{\textpartnx}[2]{{\lower1.0 ex\hbox{\begin{tikzpicture}[xscale=.3,yscale=0.3] 
\foreach \x in {1,...,#1}
{ \uvert{\x} \lvert{\x} }
#2 \end{tikzpicture}}}}
\nc{\disppartnx}[2]{{\lower1.0 ex\hbox{\begin{tikzpicture}[scale=0.3] 
\foreach \x in {1,...,#1}
{ \uvert{\x} \lvert{\x} }
#2 \end{tikzpicture}}}}
\nc{\disppartnxd}[2]{{\lower2.1 ex\hbox{\begin{tikzpicture}[scale=0.3] 
\foreach \x in {1,...,#1}
{ \uuvert{\x} \uvert{\x} \lvert{\x} }
#2 \end{tikzpicture}}}}
\nc{\disppartnxdn}[2]{{\lower2.1 ex\hbox{\begin{tikzpicture}[scale=0.3] 
\foreach \x in {1,...,#1}
{ \uuvert{\x} \lvert{\x} }
#2 \end{tikzpicture}}}}
%\nc{\disppartnxdd}[2]{{\lower2.1 ex\hbox{\begin{tikzpicture}[scale=0.3] 
%\foreach \x in {1,...,#1}
%{ \uuuvert{\x} \uuvert{\x} \uvert{\x} \lvert{\x} }
%#2 \end{tikzpicture}}}}
\nc{\disppartnxdd}[2]{{\lower3.6 ex\hbox{\begin{tikzpicture}[scale=0.3] 
\foreach \x in {1,...,#1}
{ \uuuvert{\x} \uuvert{\x} \uvert{\x} \lvert{\x} }
#2 \end{tikzpicture}}}}

\nc{\dispgax}[2]{{\lower0.0 ex\hbox{\begin{tikzpicture}[scale=0.3] 
#2
\foreach \x in {1,...,#1}
{\lvert{\x} }
 \end{tikzpicture}}}}
\nc{\textgax}[2]{{\lower0.4 ex\hbox{\begin{tikzpicture}[scale=0.3] 
#2
\foreach \x in {1,...,#1}
{\lvert{\x} }
 \end{tikzpicture}}}}
\nc{\textlinegraph}[2]{{\raise#1 ex\hbox{\begin{tikzpicture}[scale=0.8] 
#2
 \end{tikzpicture}}}}
\nc{\textlinegraphl}[2]{{\raise#1 ex\hbox{\begin{tikzpicture}[scale=0.8] 
\tikzstyle{vertex}=[circle,draw=black, fill=white, inner sep = 0.07cm]
#2
 \end{tikzpicture}}}}
\nc{\displinegraph}[1]{{\lower0.0 ex\hbox{\begin{tikzpicture}[scale=0.6] 
#1
 \end{tikzpicture}}}}
 
\nc{\disppartnthreeone}[1]{{\lower1.0 ex\hbox{\begin{tikzpicture}[scale=0.3] 
\foreach \x in {1,2,3,5,6}
{ \uvert{\x} }
\foreach \x in {1,2,4,5,6}
{ \lvert{\x} }
\draw[dotted] (3.5,1.5)--(4.5,1.5);
\draw[dotted] (2.5,0)--(3.5,0);
#1 \end{tikzpicture}}}}

\nc{\partn}[4]{\left( \begin{array}{c|c} %fine
#1 \ & \ #3 \ \ \\ \cline{2-2}
#2 \ & \ #4 \ \
\end{array} \!\!\! \right)}
\nc{\partnlong}[6]{\partn{#1}{#2}{#3,\ #4}{#5,\ #6}} %fine
\nc{\partnsh}[2]{\left( \begin{array}{c} %fine
#1 \\
#2 
\end{array} \right)}
\nc{\partncodefz}[3]{\partn{#1}{#2}{#3}{\emptyset}}
\nc{\partndefz}[3]{{\partn{#1}{#2}{\emptyset}{#3}}}
\nc{\partnlast}[2]{\left( \begin{array}{c|c}
#1 \ &  \ #2 \\
#1 \ &  \ #2
\end{array} \right)}

\nc{\uuarcx}[3]{\draw(#1,3)arc(180:270:#3) (#1+#3,3-#3)--(#2-#3,3-#3) (#2-#3,3-#3) arc(270:360:#3);}
\nc{\uuarc}[2]{\uuarcx{#1}{#2}{.4}}
\nc{\uuuarcx}[3]{\draw(#1,4.5)arc(180:270:#3) (#1+#3,4.5-#3)--(#2-#3,4.5-#3) (#2-#3,4.5-#3) arc(270:360:#3);}
\nc{\uuuarc}[2]{\uuuarcx{#1}{#2}{.4}}
\nc{\darcx}[3]{\draw(#1,0)arc(180:90:#3) (#1+#3,#3)--(#2-#3,#3) (#2-#3,#3) arc(90:0:#3);}
\nc{\darc}[2]{\darcx{#1}{#2}{.4}}
\nc{\udarcx}[3]{\draw(#1,1.5)arc(180:90:#3) (#1+#3,1.5+#3)--(#2-#3,1.5+#3) (#2-#3,1.5+#3) arc(90:0:#3);}
\nc{\udarc}[2]{\udarcx{#1}{#2}{.4}}
\nc{\uudarcx}[3]{\draw(#1,3)arc(180:90:#3) (#1+#3,3+#3)--(#2-#3,3+#3) (#2-#3,3+#3) arc(90:0:#3);}
\nc{\uudarc}[2]{\uudarcx{#1}{#2}{.4}}
\nc{\uarcx}[3]{\draw(#1,1.5)arc(180:270:#3) (#1+#3,1.5-#3)--(#2-#3,1.5-#3) (#2-#3,1.5-#3) arc(270:360:#3);}
\nc{\uarc}[2]{\uarcx{#1}{#2}{.4}}
\nc{\darcxhalf}[3]{\draw(#1,0)arc(180:90:#3) (#1+#3,#3)--(#2,#3) ;}
\nc{\darchalf}[2]{\darcxhalf{#1}{#2}{.4}}
\nc{\uarcxhalf}[3]{\draw(#1,1.5)arc(180:270:#3) (#1+#3,1.5-#3)--(#2,1.5-#3) ;}
\nc{\uarchalf}[2]{\uarcxhalf{#1}{#2}{.4}}
\nc{\uarcxhalfr}[3]{\draw (#1+#3,1.5-#3)--(#2-#3,1.5-#3) (#2-#3,1.5-#3) arc(270:360:#3);}
\nc{\uarchalfr}[2]{\uarcxhalfr{#1}{#2}{.4}}

\nc{\bdarcx}[3]{\draw[blue](#1,0)arc(180:90:#3) (#1+#3,#3)--(#2-#3,#3) (#2-#3,#3) arc(90:0:#3);}
\nc{\bdarc}[2]{\darcx{#1}{#2}{.4}}
\nc{\rduarcx}[3]{\draw[red](#1,0)arc(180:270:#3) (#1+#3,0-#3)--(#2-#3,0-#3) (#2-#3,0-#3) arc(270:360:#3);}
\nc{\rduarc}[2]{\uarcx{#1}{#2}{.4}}
\nc{\duarcx}[3]{\draw(#1,0)arc(180:270:#3) (#1+#3,0-#3)--(#2-#3,0-#3) (#2-#3,0-#3) arc(270:360:#3);}
\nc{\duarc}[2]{\uarcx{#1}{#2}{.4}}

\nc{\uv}[1]{\fill (#1,2)circle(.1);}
\nc{\lv}[1]{\fill (#1,0)circle(.1);}
\nc{\stline}[2]{\draw(#1,2)--(#2,0);}
\nc{\tlab}[2]{\draw(#1,2)node[above]{\tiny $#2$};}
\nc{\tudots}[1]{\draw(#1,2)node{$\cdots$};}
\nc{\tldots}[1]{\draw(#1,0)node{$\cdots$};}

\nc{\uvw}[1]{\fill[white] (#1,2)circle(.1);}
\nc{\huv}[1]{\fill (#1,1)circle(.1);}
\nc{\llv}[1]{\fill (#1,-2)circle(.1);}
\nc{\arcup}[2]{
\draw(#1,2)arc(180:270:.4) (#1+.4,1.6)--(#2-.4,1.6) (#2-.4,1.6) arc(270:360:.4);
}
\nc{\harcup}[2]{
\draw(#1,1)arc(180:270:.4) (#1+.4,.6)--(#2-.4,.6) (#2-.4,.6) arc(270:360:.4);
}
\nc{\arcdn}[2]{
\draw(#1,0)arc(180:90:.4) (#1+.4,.4)--(#2-.4,.4) (#2-.4,.4) arc(90:0:.4);
}
\nc{\cve}[2]{
\draw(#1,2) to [out=270,in=90] (#2,0);
}
\nc{\hcve}[2]{
\draw(#1,1) to [out=270,in=90] (#2,0);
}
\nc{\catarc}[3]{
\draw(#1,2)arc(180:270:#3) (#1+#3,2-#3)--(#2-#3,2-#3) (#2-#3,2-#3) arc(270:360:#3);
}

\nc{\arcr}[2]{
\draw[red](#1,0)arc(180:90:.4) (#1+.4,.4)--(#2-.4,.4) (#2-.4,.4) arc(90:0:.4);
}
\nc{\arcb}[2]{
\draw[blue](#1,2-2)arc(180:270:.4) (#1+.4,1.6-2)--(#2-.4,1.6-2) (#2-.4,1.6-2) arc(270:360:.4);
}
\nc{\loopr}[1]{
\draw[blue](#1,-2) edge [out=130,in=50,loop] ();
}
\nc{\loopb}[1]{
\draw[red](#1,-2) edge [out=180+130,in=180+50,loop] ();
}
%\nc{\arcr}[2]{
%\draw[red](#1,0-2)arc(180:90:.4) (#1+.4,.4-2)--(#2-.4,.4-2) (#2-.4,.4-2) arc(90:0:.4);
%}
%\nc{\arcb}[2]{
%\draw[blue](#1,2-2-2)arc(180:270:.4) (#1+.4,1.6-2-2)--(#2-.4,1.6-2-2) (#2-.4,1.6-2-2) arc(270:360:.4);
%}
%\nc{\loopr}[1]{
%\draw[red](#1,0-2) edge [out=130,in=50,loop] ();
%}
%\nc{\loopb}[1]{
%\draw[blue](#1,0-2) edge [out=180+130,in=180+50,loop] ();
%}
\nc{\redto}[2]{\draw[red](#1,0)--(#2,0);}
\nc{\bluto}[2]{\draw[blue](#1,0)--(#2,0);}
\nc{\dotto}[2]{\draw[dotted](#1,0)--(#2,0);}
\nc{\lloopr}[2]{\draw[red](#1,0)arc(0:360:#2);}
\nc{\lloopb}[2]{\draw[blue](#1,0)arc(0:360:#2);}
\nc{\rloopr}[2]{\draw[red](#1,0)arc(-180:180:#2);}
\nc{\rloopb}[2]{\draw[blue](#1,0)arc(-180:180:#2);}
\nc{\uloopr}[2]{\draw[red](#1,0)arc(-270:270:#2);}
\nc{\uloopb}[2]{\draw[blue](#1,0)arc(-270:270:#2);}
\nc{\dloopr}[2]{\draw[red](#1,0)arc(-90:270:#2);}
\nc{\dloopb}[2]{\draw[blue](#1,0)arc(-90:270:#2);}
\nc{\llloopr}[2]{\draw[red](#1,0-2)arc(0:360:#2);}
\nc{\llloopb}[2]{\draw[blue](#1,0-2)arc(0:360:#2);}
\nc{\lrloopr}[2]{\draw[red](#1,0-2)arc(-180:180:#2);}
\nc{\lrloopb}[2]{\draw[blue](#1,0-2)arc(-180:180:#2);}
\nc{\ldloopr}[2]{\draw[red](#1,0-2)arc(-270:270:#2);}
\nc{\ldloopb}[2]{\draw[blue](#1,0-2)arc(-270:270:#2);}
\nc{\luloopr}[2]{\draw[red](#1,0-2)arc(-90:270:#2);}
\nc{\luloopb}[2]{\draw[blue](#1,0-2)arc(-90:270:#2);}

\nc{\larcb}[2]{
\draw[blue](#1,0-2)arc(180:90:.4) (#1+.4,.4-2)--(#2-.4,.4-2) (#2-.4,.4-2) arc(90:0:.4);
}
\nc{\larcr}[2]{
\draw[red](#1,2-2-2)arc(180:270:.4) (#1+.4,1.6-2-2)--(#2-.4,1.6-2-2) (#2-.4,1.6-2-2) arc(270:360:.4);
}

\rnc{\H}{\mathscr H}
\rnc{\L}{\mathscr L}
\nc{\R}{\mathscr R}
\nc{\D}{\mathcal D}
\nc{\J}{\mathscr D}

\nc{\ssim}{\mathrel{\raise0.25 ex\hbox{\oalign{$\approx$\crcr\noalign{\kern-0.84 ex}$\approx$}}}}
\nc{\POI}{\mathcal{POI}}
\nc{\wb}{\overline{w}}
\nc{\ub}{\overline{u}}
\nc{\vb}{\overline{v}}
\nc{\fb}{\overline{f}}
\nc{\gb}{\overline{g}}
\nc{\hb}{\overline{h}}
\nc{\pb}{\overline{p}}
\rnc{\sb}{\overline{s}}
\nc{\XR}{\pres{X}{R\,}}
\nc{\YQ}{\pres{Y}{Q}}
\nc{\ZP}{\pres{Z}{P\,}}
\nc{\XRone}{\pres{X_1}{R_1}}
\nc{\XRtwo}{\pres{X_2}{R_2}}
\nc{\XRthree}{\pres{X_1\cup X_2}{R_1\cup R_2\cup R_3}}
\nc{\er}{\eqref}
\nc{\larr}{\mathrel{\hspace{-0.35 ex}>\hspace{-1.1ex}-}\hspace{-0.35 ex}}
\nc{\rarr}{\mathrel{\hspace{-0.35 ex}-\hspace{-0.5ex}-\hspace{-2.3ex}>\hspace{-0.35 ex}}}
\nc{\lrarr}{\mathrel{\hspace{-0.35 ex}>\hspace{-1.1ex}-\hspace{-0.5ex}-\hspace{-2.3ex}>\hspace{-0.35 ex}}}
\nc{\nn}{\tag*{}}
\nc{\epfal}{\tag*{$\Box$}}
\nc{\tagd}[1]{\tag*{(#1)$'$}}
\nc{\ldb}{[\![}
\nc{\rdb}{]\!]}
\nc{\sm}{\setminus}
\nc{\I}{\mathcal I}
\nc{\InSn}{\I_n\setminus\S_n}
%\nc{\dom}{\operatorname{dom}_{\operatorname{u}}} \nc{\codom}{\operatorname{dom}_{\operatorname{l}}}
%\nc{\dom}{\operatorname{dom}_u} \nc{\codom}{\operatorname{dom}_l}
\nc{\dom}{\operatorname{dom}^\wedge} \nc{\codom}{\operatorname{dom}_\vee}%\nc{\codom}{\operatorname{im}}
\nc{\ojin}{1\leq j<i\leq n}
%\nc{\R}{\mathcal R}
%\rnc{\L}{\mathcal L}
\nc{\eh}{\widehat{e}}
\nc{\wh}{\widehat{w}}
\nc{\uh}{\widehat{u}}
\nc{\vh}{\widehat{v}}
\nc{\sh}{\widehat{s}}
\nc{\fh}{\widehat{f}}
\nc{\textres}[1]{\text{\emph{#1}}}
\nc{\aand}{\emph{\ and \ }}
\nc{\iif}{\emph{\ if \ }}
\nc{\textlarr}{\ \larr\ }
\nc{\textrarr}{\ \rarr\ }
\nc{\textlrarr}{\ \lrarr\ }

\nc{\comma}{,\ }

\nc{\COMMA}{,\quad}
\nc{\TnSn}{\T_n\setminus\S_n} 
\nc{\TmSm}{\T_m\setminus\S_m} 
\nc{\TXSX}{\T_X\setminus\S_X} 
\rnc{\S}{\mathcal S}

\nc{\T}{\mathcal T} 
\nc{\A}{\mathscr A} 
\nc{\B}{\mathcal B} 
\rnc{\P}{\mathcal P} 
\nc{\K}{\mathcal K}
\nc{\PB}{\mathcal{PB}} 
\nc{\rank}{\operatorname{rank}}

\nc{\mtt}{\!\!\!\mt\!\!\!}

\nc{\sub}{\subseteq}
\nc{\la}{\langle}
\nc{\ra}{\rangle}
\nc{\mt}{\mapsto}
\nc{\im}{\mathrm{im}}
\nc{\id}{\mathrm{id}}
\nc{\bn}{\mathbf{n}}
\nc{\ba}{\mathbf{a}}
\nc{\bl}{\mathbf{l}}
\nc{\bm}{\mathbf{m}}
\nc{\bk}{\mathbf{k}}
\nc{\br}{\mathbf{r}}
\nc{\ve}{\varepsilon}
\nc{\al}{\alpha}
\nc{\be}{\beta}
\nc{\ga}{\gamma}
\nc{\Ga}{\Gamma}
\nc{\de}{\delta}
\nc{\ka}{\kappa}
\nc{\lam}{\lambda}
\nc{\Lam}{\Lambda}
\nc{\si}{\sigma}
\nc{\Si}{\Sigma}
\nc{\oijn}{1\leq i<j\leq n}
\nc{\oijm}{1\leq i<j\leq m}

\nc{\comm}{\rightleftharpoons}
\nc{\AND}{\qquad\text{and}\qquad}

\nc{\bit}{\vspace{-3 truemm}\begin{itemize}}
\nc{\bmc}{\vspace{-3 truemm}\begin{multicols}}
\nc{\emc}{\end{multicols}\vspace{-3 truemm}}
\nc{\eit}{\end{itemize}\vspace{-3 truemm}}
\nc{\ben}{\vspace{-3 truemm}\begin{enumerate}}
\nc{\een}{\end{enumerate}\vspace{-3 truemm}}
\nc{\eitres}{\end{itemize}}

\nc{\set}[2]{\{ {#1} : {#2} \}} 
\nc{\bigset}[2]{\big\{ {#1}: {#2} \big\}} 
\nc{\Bigset}[2]{\Big\{ \,{#1}\, \,\Big|\, \,{#2}\, \Big\}}

\nc{\pres}[2]{\la {#1} \,|\, {#2} \ra}
\nc{\bigpres}[2]{\big\la {#1} \,\big|\, {#2} \big\ra}
\nc{\Bigpres}[2]{\Big\la \,{#1}\, \,\Big|\, \,{#2}\, \Big\ra}
\nc{\Biggpres}[2]{\Bigg\la {#1} \,\Bigg|\, {#2} \Bigg\ra}

\nc{\pf}{\noindent{\bf Proof.}  }
\nc{\epf}{\hfill$\Box$\bigskip}
\nc{\epfres}{\hfill$\Box$}
\nc{\pfnb}{\pf}
\nc{\epfnb}{\bigskip}
\nc{\pfthm}[1]{\bigskip \noindent{\bf Proof of Theorem \ref{#1}}\,\,  } 
\nc{\pfprop}[1]{\bigskip \noindent{\bf Proof of Proposition \ref{#1}}\,\,  } 
%\nc{\pfthm}{\noindent{\bf Proof of Theorem \ref{mainthm} modulo Propositions \ref{prop1} and \ref{prop2}}\,\,  } 
\nc{\epfreseq}{\tag*{$\Box$}}

\makeatletter
\newcommand\footnoteref[1]{\protected@xdef\@thefnmark{\ref{#1}}\@footnotemark}
\makeatother

\numberwithin{equation}{section}

\newtheorem{thm}[equation]{Theorem}
\newtheorem{lemma}[equation]{Lemma}
\newtheorem{cor}[equation]{Corollary}
\newtheorem{prop}[equation]{Proposition}

\theoremstyle{definition}

\newtheorem{rem}[equation]{Remark}
\newtheorem{defn}[equation]{Definition}
\newtheorem{eg}[equation]{Example}
\newtheorem{ass}[equation]{Assumption}

%\title{The singular part of the endomorphism monoid of a finite partition~\vspace{-5ex}}
\title{Idempotent generation in the endomorphism monoid of a uniform partition} %~\vspace{-5ex}}
\author{
Igor Dolinka%
\footnote{The first named author gratefully acknowledges the support of Grant No.~174019 of the Ministry of Education, Science, and Technological Development of the Republic of Serbia, and Grant No.~1136/2014 of the Secretariat of Science and Technological Development of the Autonomous Province of Vojvodina.}
\\
{\footnotesize \emph{Department of Mathematics and Informatics}}\\
{\footnotesize \emph{University of Novi Sad, Trg Dositeja Obradovi\'ca 4, 21101 Novi Sad, Serbia}}\\
{\footnotesize {\tt dockie@dmi.uns.ac.rs}}\\~\\
James East\\
{\footnotesize \emph{Centre for Research in Mathematics; School of Computing, Engineering and Mathematics}}\\
{\footnotesize \emph{University of Western Sydney, Locked Bag 1797, Penrith NSW 2751, Australia}}\\
{\footnotesize {\tt J.East\,@\,uws.edu.au}}
}

\maketitle
%\begin{center}
%{\large Igor Dolinka\footnote{Department of Mathematics and Informatics, University of Novi Sad, Trg Dositeja Obradovi\'ca 4, 21101 Novi Sad, Serbia. {\it Email:} {\tt dockie@dmi.uns.ac.rs}} and James East\footnote{Centre for Research in Mathematics, School of Computing, Engineering and Mathematics, University of Western Sydney, Locked Bag 1797, Penrith NSW 2751, Australia. {\it Email:} {\tt J.East@uws.edu.au}}}\\
%~\\ 
%{\large \today}
%\end{center}

~\vspace{-10ex}
\begin{abstract}
Denote by $\mathcal T_n$ and $\mathcal S_n$ the full transformation semigroup and the symmetric group on the set $\{1,\ldots,n\}$, and $\mathcal E_n=\{1\}\cup(\mathcal T_n\setminus\mathcal S_n)$.  Let $\mathcal T(X,\mathcal P)$ denote the monoid of all transformations of the finite set $X$ preserving a uniform partition $\mathcal P$ of $X$ into $m$ subsets of size $n$, where $m,n\geq2$.  We enumerate the idempotents of $\mathcal T(X,\mathcal P)$, and describe the submonoid $S=\langle E\rangle$ generated by the idempotents $E=E(\mathcal T(X,\mathcal P))$.  We show that $S=S_1\cup S_2$, where $S_1$ is a direct product of $m$ copies of $\mathcal E_n$, and $S_2$ is a wreath product of $\mathcal T_n$ with $\mathcal T_m\setminus\mathcal S_m$.  We calculate the rank and idempotent rank of $S$, showing that these are equal, and we also classify and enumerate all the idempotent generating sets of minimal size.  In doing so, we also obtain new results about arbitrary idempotent generating sets of~$\mathcal E_n$.


\end{abstract}

\section{Introduction}\label{sect:intro}

Let $M$ be a monoid and $E(M)=\set{x\in M}{x^2=x}$ the set of all idempotents of $M$.  For a subset $U\sub M$, we write $\la U\ra$ (respectively, $\la U\rasgp$) for the submonoid (respectively, subsemigroup) of $M$ generated by $U$, which consists of all products $x_1\cdots x_k$ where $k\geq0$ (respectively, $k\geq1$) and $x_1,\ldots,x_k\in U$.  (By convention, the empty product is equal to the identity element $1\in M$.)  The \emph{rank} of $M$ is the minimal cardinality of a subset $U\sub M$ such that $M=\la U\ra$.  If $M$ is idempotent generated, then the \emph{idempotent rank} of $M$, denoted $\idrank(M)$, is the minimal cardinality of a subset $U\sub E(M)$ such that $M=\la U\ra$.  All monoids we consider will have an irreducible identity element $1$;  in other words, the only solution to the equation $xy=1$ is $x=y=1$.  For such a monoid, the smallest semigroup generating set has size $1+\rank(M)$, with a similar statement holding for idempotent generating sets.  %minimal cardinality of a subset $U\sub M$ such that $M=\la U\rasgp$ is equal to $1+\rank(M)$.

%For such a monoid $M$, we see that $S=M\sm\{1\}$ is a subsemigroup of $M$, and that $M=\la U\ra^1$ if and only if $S=\la U\ra$.  

The \emph{full transformation semigroup} on a set $X$, denoted $\T_X$, is the set of all transformations of~$X$ (i.e., all functions $X\to X$), under the semigroup operation of composition.  The group of units of $\T_X$ is the \emph{symmetric group} $\S_X$, consisting of all permutations of $X$ (i.e., all bijections $X\to X$).  We also write $\E_X=\la E(\T_X)\ra$ for the idempotent generated submonoid of $\T_X$.  
When $X=\bn=\{1,\ldots,n\}$, we write $\T_X=\T_n$, and similarly for $\S_n$ and $\E_n$.  (Note that $\T_0=\S_0=\E_0$ has a single element; namely, the empty map~$\emptyset$.)  A much celebrated result of Howie \cite{Howie66} states that ${\E_n=\{1\}\cup(\TnSn)}$, where $1\in\T_n$ denotes the identity mapping.  (Infinite~$\E_X$ was also described in \cite{Howie66}.)  In subsequent work, Howie \cite{Howie1978} showed that $\idrank(\Sing_n)={n\choose2}$ if $n\geq3$%
%(while $\idrank(\Sing_2)=2$, and $\idrank(\Sing_n)=0$ if $n\leq1$)
, and classified the minimal idempotent generating sets of $\Sing_n$, showing that these are in one-one correspondence with the strongly connected tournaments on $n$ vertices; the enumeration of such tournaments was given by Wright \cite{Wright1970}.  Gomes and Howie \cite{Gomes1987} showed that also $\rank(\E_n)=\idrank(\E_n)$ for all $n$.  These results initiated a vibrant direction for research in combinatorial semigroup theory.  For example, they have been extended to semigroups of matrices \cite{Erdos1967, laffey83,Gray2007,Fountain1991}, endomorphisms of (finite and infinite dimensional) independence algebras \cite{Fountain1992,Fountain1993}, and partitions \cite{Maltcev2007,East2011_2,EF,EastGray}.

Now let $X$ be an arbitrary set and $\P=\set{C_i}{i\in I}$ a partition of $X$; that is, the sets $C_i$ are non-empty, pairwise disjoint, and their union is all of $X$.  The set
\[
\TXP = \set{f\in\T_X}{(\forall i\in I)(\exists j\in I) \ C_if\sub C_j},
\]
consisting of all transformations of $X$ preserving $\P$,
is a submonoid of $\T_X$, and may be thought of as the set of all continuous mappings with respect to the topology with basis $\P$.  This monoid appears to have been introduced by Huisheng Pei \cite{Pei1994}, who later characterised the regular elements of $\TXP$ and described Green's relations \cite{Pei2005}; see also \cite{Pei2009}.  In \cite{Pei2005b}, Pei raised the question of determining $\rank(\TXP)$ in the case that the set $X$ is finite and the partition~$\P$ is uniform (i.e., $|C_i|=|C_j|$ for each $i,j$), and obtained an upper bound of $6$ for this rank.  The question was settled in \cite{AS2009} where it was shown that $\rank(\TXP)=4$, making use of results on wreath products of symmetric groups and transformation semigroups; see also \cite{ABMS2014} for the calculation of $\rank(\TXP)$ for an arbitrary partition $\P$ of a finite set $X$.

The purpose of the current work is to consider the idempotent generated subsemigroup $$\EXP=\la E(\TXP)\ra$$ of $\TXP$ where $X$ is finite and $\P$ is uniform.  Our results include:
\ben[(i)]
\item characterisation and enumeration of the idempotents of $\TXP$ --- see Propositions~\ref{idempotent_prop} and \ref{etxp_rec},
\item description of the elements and structure of $\EXP$ --- see Proposition \ref{prop:S1S2}, 
\item calculation of $\rank(\EXP)$ and $\idrank(\EXP)$ --- see Theorem \ref{rank_thm},
\item classification and enumeration of the minimal idempotent generating sets of $\EXP$ --- see Theorem \ref{enum_thm}.
\een
Note that by a ``minimal generating set'', we always mean a generating set of minimal size, rather than a generating set that contains no smaller generating set, but we will show that these two notions of minimality are equivalent in the context of idempotent generating sets of $\EXP$ --- see Theorem \ref{enum_thm}.  In order to obtain the results alluded to in (iv), it is necessary to extend Howie's results from \cite{Howie1978} and enumerate the idempotent generating sets of $\E_n$ that are not necessarily minimal; we believe these results are interesting in their own right --- see Theorem \ref{wnk}.  Using the current article as a starting point, many of the above-mentioned results are extended to the non-uniform case in \cite{DEM2015}.

\section{Preliminaries on $\T_n$}\label{sect:TX}

In this section, we record the results concerning $\T_n$, $E(\T_n)$ and $\E_n=\la E(\T_n)\ra$ that we will need in what follows.  %Most are well known, but we include a couple of new results that we think are interesting enough in their own right.

%The group of units of $\T_X$ is the symmetric group $\S_X$, which consists of all permutations of $X$.  We write $\Sing(X)$ for the set $(\TXSX)\cup\{1\}$, and we call this the \emph{singular submonoid} of $\T_X$.  (Note that it is customary not to include the identity transformation $1\in\T_X$ in $\Sing(X)$, but it is convenient to do so here.)

Recall that the \emph{rank} of a transformation $f\in\T_X$ is $\rank(f)=|\im(f)|$, the cardinality of the image $\im(f)$ of $f$.  For $k\in\bn$, we will write $D_{nk}=\set{f\in\T_n}{\rank(f)=k}$.\footnote{These sets form the so-called \emph{$\mathscr D$-classes} of $\T_n$%
%, and they form a chain, $D_{n1}<D_{n2}<\cdots<D_{nn}$
.  No knowledge of Green's relations, which include the $\mathscr D$ relation, will be assumed but the reader may refer to a monograph such as \cite{Howie,Hig} for details if they wish.}  For a subset $U$ of a semigroup $S$, we write $E(U)=U\cap E(S)$ for the set of all idempotents in~$U$.  The next result is folklore and is easily checked.

\begin{prop}\label{ETn_prop}
A transformation $f\in\T_X$ is an idempotent if and only if $f$ acts as the identity on its image.  If $n\geq1$, then 
\ben
\itemit{i} $|E(D_{nk})|={n\choose k}k^{n-k}$ for $1\leq k\leq n$, and
\itemit{ii} $|E(\T_n)| = \sum_{k=1}^n{n\choose k}k^{n-k}$. \epfres
\een
\end{prop}

\begin{figure}[ht]
\begin{center}
\begin{tikzpicture}[scale=0.6]
\uv1
\uv3
\uv4
\uv5
\uv7
\uv8
\uv9
\uv{11}
\lv1
\lv3
\lv4
\lv5
\lv7
\lv8
\lv9
\lv{11}
\stline11
\stline33
\stline44
\stline55
\stline77
\stline84
\stline99
\stline{11}{11}
\tlab11
\tlab4i
\tlab8j
\tlab{11}n
\tudots2
\tudots6
\tudots{10}
\tldots2
\tldots6
\tldots{10}
\end{tikzpicture}
\qquad\qquad
\begin{tikzpicture}[scale=0.6]
\uv1
\uv3
\uv4
\uv5
\uv7
\uv8
\uv9
\uv{11}
\lv1
\lv3
\lv4
\lv5
\lv7
\lv8
\lv9
\lv{11}
\stline11
\stline33
\stline48
\stline55
\stline77
\stline88
\stline99
\stline{11}{11}
\tlab11
\tlab4j
\tlab8i
\tlab{11}n
\tudots2
\tudots6
\tudots{10}
\tldots2
\tldots6
\tldots{10}
\end{tikzpicture}
    \caption{The idempotent $e_{ij}\in E(\D_n)$ in the cases $i<j$ (left) and $i>j$ (right).}
    \label{fig:eij}
   \end{center}
 \end{figure}
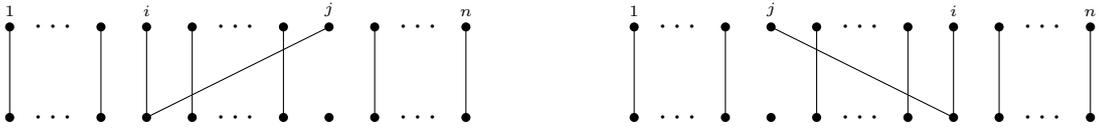

From now on, we will write $\D_n=D_{n,n-1}$ for $n\geq2$.  By convention, we also define $\D_1=\D_0=\emptyset$.  By the previous proposition, 
\[
E(\D_n) = \set{e_{ij}}{i,j\in \bn,\ i\not=j},
\]
where $e_{ij}\in\T_n$ denotes the transformation that maps $j$ to $i$ and maps the rest of $\bn$ identically; see Figure \ref{fig:eij} for an illustration.  

\begin{thm}[Howie, \cite{Howie66}]
If $n\geq0$, then $  \Sing_n=\la E(\D_n)\ra =\{1\}\cup(\TnSn)$.  \epf
%If $n\geq0$, then $\la E(\T_n)\ra = \la E(\D_n)\ra = \Sing_n=\{1\}\cup(\TnSn)$.  \epf
%, and 
%\[
%\rank(\Sing_n)=\idrank(\Sing_n)=\begin{cases} 2 &\text{if $n=2$} \\ 
%{n\choose2} &\text{if $n\geq3$.} \end{cases}
%\]
\end{thm}

\begin{rem}\label{min_rem}
Note that the empty map $\emptyset\in\T_0$ is the identity element of $\T_0$.  Since the set $\set{f\in\T_n}{\rank(f)\leq n-2}$ is an ideal of $\T_n$, and since $\E_n=\la\D_n\ra$, it follows that whenever $U$ is a generating set for $\E_n$, then so too is $U\cap\D_n$.  In particular, any minimal generating set for $\E_n$ is contained in $\D_n$, and any minimal idempotent generating set is contained in $E(\D_n)$.  See \cite{JEtnsn2} for a presentation for $\E_n$ with respect to the generating set $E(\D_n)$.
\end{rem}

The rank and idempotent rank of $\E_n$ are calculated in \cite{Gomes1987,Howie1978}, and the (minimal) idempotent generating sets are classified in \cite{Howie1978}.
%In \cite{Gomes1987,Howie1978}, Howie calculated the rank and idempotent rank of $\E_n$, and classified the (minimal) idempotent generating sets.  
(A classification of all the generating sets of $\E_n$ may be found in \cite{Ayik2013}.)  We now describe Howie's classification, as it features prominently in our results in subsequent sections.

%Consider a subset $U\sub E(\D_n)$.  Define a digraph $\Ga_U$ with vertex set $\bn$ and a directed edge $i\to j$ for each $e_{ij}\in U$.  
We say that a digraph $\Ga$ on vertex set $V$ is \emph{complete} if its underlying undirected graph (obtained by changing each directed edge $u\to v$ to an undirected edge $u-v$ and replacing any double edges with single edges) is the complete graph on $V$, and we say $\Ga$ is \emph{strongly connected} if $|V|=1$ or, for any $u,v\in V$, there is a directed path from $u$ to $v$ in $\Ga$.  In what follows, we interpret ${i\choose j}=0$ if $j>i$.

%We say that $\Ga_F$ is \emph{complete} if its underlying un-directed graph is the complete graph on $\bn$, and we say $\Ga_F$ is \emph{strongly connected} if, for any $x,y\in X$, there is a directed path from $x$ to $y$ in $\Ga_F$.

\begin{thm}[{Howie, \cite{Howie1978}}]\label{Howie1978_thm}
Let $U\sub E(\D_n)$ where $n\geq0$, and define a digraph $\Ga_U$ with vertex set $\bn$ and a directed edge $i\to j$ for each $e_{ij}\in U$.  Then $\E_n=\la U\ra$ if and only if the graph $\Ga_U$ is strongly connected and complete.  Further,
%If $n\geq2$, then $\la E(\T_n)\ra = \Sing_n=\{1\}\cup(\TnSn)$, and 
\[\epfreseq
\rank(\Sing_n)=\idrank(\Sing_n)=\begin{cases} 2 &\text{if $n=2$} \\ 
{n\choose2} &\text{otherwise.} \end{cases}
%{n\choose2} &\text{if $n\geq3$.} \end{cases}
\]
\end{thm}

Note that when $n\leq1$, $\E_n=\{1\}$ so that, indeed, $\rank(\E_n)=\idrank(\E_n)=0={n\choose2}$.  There is a unique strongly connected, complete graph on vertex set ${\bf 2}=\{1,2\}$, so there is a unique (minimal idempotent) generating set of $\E_2$, as may also be easily verified directly.  

A digraph $\Ga$ on vertex set $V$ is called a \emph{tournament} if, for each $u,v\in V$ with $u\not=v$, $\Ga$ contains exactly one of the edges $u\to v$ or $v\to u$.  Part (i) of the next result is due to Howie \cite{Howie1978}, and part (ii) to Wright \cite{Wright1970}.

%The semigroup $\E_2=\{e_{12},e_{21}\}$ is a right zero semigroup, and hence it has a unique (idempotent) generating set; namely, itself.  But for $n\geq3$, 

\begin{thm}[Howie \cite{Howie1978} and Wright \cite{Wright1970}]\label{wn}
\leavevmode\newline\vspace{-5mm}
\bit
\itemit{i} The minimal idempotent generating sets of $\E_n$ with $n\not=2$ are in one-one correspondence with the strongly connected tournaments on $n$ vertices.  (There is a unique minimal idempotent generating set of $\E_2$.)
\itemit{ii} Let $w_n$ denote the number of strongly connected tournaments on $n\geq0$ vertices.  Then
\[
w_0=1,\qquad w_n=F_n - \sum_{s=1}^{n-1}{n\choose s}w_{s}F_{n-s} \quad\text{for $n\geq1$},
\]
where $F_n=2^{{n\choose 2}}=2^{n(n-1)/2}$.  \epfres
\eit

\end{thm}

In the course of our investigations, we will also need to know the total number of generating sets of $\Sing_n$ consisting of $k$ idempotents from $\D_n$, where $k$ is not necessarily equal to the minimal size of~${n\choose2}$.  

Let $\Ga$ be a complete digraph on vertex set $V$.  We say that $\Ga$ has a double edge $u-v$ if both $u\to v$ and $v\to u$ are edges of $\Ga$.  Define a relation $\sim_\Ga$ on the vertex set $V$ by $u\sim_\Ga v$ if $u=v$ or there is a directed path from $u$ to $v$ in $\Ga$ and one from $v$ to $u$.  The equivalence classes with respect to $\sim_\Ga$ are the \emph{strongly connected components} of $\Ga$.  The set $V/\!\!\sim_\Ga$ of strongly connected components is totally ordered: if $A,B$ are two strongly connected components, we say $A>B$ if every edge between a vertex from $A$ and a vertex from $B$ points from the vertex from $A$ to the vertex from $B$.  Because this is a total order, there is a maximal (and minimal) strongly connected component.

\begin{thm}\label{wnk}
For $n\geq0$ and $0\leq k\leq {n\choose 2}$, let $w_{nk}$ denote the number of strongly connected, complete digraphs on vertex set $\bn$ that have $k$ double edges.  
\bit
\itemit{i} For $0\leq k\leq{n\choose2}$, $w_{nk}$ is equal to the number of subsets $U\sub E(\D_n)$ such that $\la U\ra=\E_n$ and $|U|={n\choose2}+k$.
\vspace{-5mm}
%generating sets of $\Sing_n$ from $E(\D_n)$ of size ${n\choose2}+k$ is equal to $w_{nk}$.
\itemit{ii} The number of subsets $U\sub E(\D_n)$ such that $\la U\ra=\Sing_n$ is equal to $\displaystyle{\sum_{k=0}^{{n\choose2}}w_{nk}}$.
%The number of generating sets of $\Sing_n$ from $E(\D_n)$ is equal to $\displaystyle{\sum_{k=0}^{{n\choose2}}w_{nk}}$.
\eit
The numbers $w_{nk}$ satisfy the recurrence 
\[
w_{00}=1,\qquad w_{nk} = F_{nk} - \sum_{s=1}^{n-1}{n\choose s} \sum_{l=0}^{k} w_{sl}F_{n-s,k-l} \quad\text{for $n\geq1$},
%w_{20}=0,\ w_{21}=1,\quad w_{nk} = F_{nk} - \sum_{s=1}^{n-1}{n\choose s} \sum_{l=0}^{k} w_{sl}F_{n-s,k-l},
\]
where $F_{nk} = \displaystyle{{{n\choose2}\choose k}\cdot2^{{n\choose 2}-k}}$.
\end{thm}

\pf Parts (i) and (ii) are clear, based on Theorem \ref{Howie1978_thm}, so it suffices to prove the recurrence for~$w_{nk}$.  The value for $n=0$ is clear.  Now suppose $n\geq1$.  Note that $F_{nk}$ is the total number of complete digraphs on vertex set $\bn$ with $k$ double edges, since to specify such a graph, we must choose the double edges in ${{n\choose2}\choose k}$ ways, and then choose the orientation of the remaining edges in $2^{{n\choose2}-k}$ ways.  From this value, we subtract the number of complete digraphs on vertex set $\bn$ with $k$ double edges that are \emph{not} strongly connected.  Clearly, there are no such graphs if $n=1$ (in which case the stated recurrence says $w_{10}=F_{10}=1$), so we assume $n\geq2$ for the remainder of the proof.  Consider a complete, but not strongly connected, graph $\Ga$ on vertex set~$\bn$ with $k$ double edges.  To specify $\Ga$, we first choose a subset $\emptyset\subsetneq M\subsetneq\bn$ to be the maximum strongly connected component in the total order on $\bn/\!\!\sim_\Ga$; this may be done in ${n\choose s}$ ways, where $|M|=s$ (with $1\leq s\leq n-1$).  We then choose the edges within $M$ in $w_{sl}$ ways for some $0\leq l\leq k$, the edges within $\bn\sm M$ in $F_{n-s,k-l}$ ways, and the remaining edges all point from a vertex in $M$ to a vertex in $\bn\sm M$.  Summing over all $s,l$ gives the desired result. \epf

\begin{rem}\label{idempotent_rem}
Note that $F_{nk}=0$ if $k>{n\choose2}$.  We will not need part (ii) of the previous theorem, but it is included for completeness.  Since any idempotent generating set for $\E_n$ contains a generating set consisting of idempotents from $\D_n$ (as noted in Remark \ref{min_rem}), it follows that the total number of idempotent generating sets of $\E_n$ is equal to 
\[
\sum_{k=0}^{{n\choose2}}w_{nk} \times 2^{1+\sum_{l=1}^{n-2}{n\choose l}l^{n-l}}.
\]
(The above expression concerns \emph{monoid} idempotent generating sets; for the total number of \emph{semigroup} idempotent generating sets, we must divide by $2$.)
Note also that $w_{n0}=w_n$ and $F_{n0}=F_n$, so the $k=0$ case of the previous result is Theorem \ref{wn}.  Calculated values of $w_{nk}$ and $\sum_{k=0}^{{n\choose2}}w_{nk}$ may be found in Tables \ref{tab:wnk} and \ref{tab:wnksum}.  Table \ref{tab:wnk} also includes the values of $w_n=w_{n0}$ in the first column.
\end{rem}

\begin{table}[h]
\begin{center}
{
\begin{tabular}{|c||rrrrrrrrrrr|}
\hline
$n\sm k$ & 0 & 1 & 2 & 3 & 4 & 5 & 6 & 7 & 8 & 9 & 10 \\ 
\hline\hline
0     &      1&&&&&&&&&&\\
1     &      1&&&&&&&&&&\\
2     &      0&1&&&&&&&&&\\
3     &      2&6&6&1&&&&&&&\\
4     &      24&108&186&152&60&12&1&&&&\\
5     &     544&3400&9090&13660&12820&7944&3350&960&180&20&1\\
\hline
\end{tabular}
}
\end{center}
\caption{Calculated values of $w_{nk}$; see Theorem \ref{wnk}(i) for more details.
% which gives the total number of idempotent generating sets of $\E_n$ from $E(\D_n)$ of size ${n\choose2}+k$.
}
\label{tab:wnk}
\end{table}

\begin{table}[h]
\begin{center}
{
\begin{tabular}{|c||rrrrrrrrr|}
\hline
$n$ & 0 & 1 & 2 & 3 & 4 & 5 & 6 & 7  & 8 \\ 
\hline\hline
     &      1&1&1&15&543&51969&13639329&10259025615&22709334063807\\
\hline
\end{tabular}
}
\end{center}
\caption{Calculated values of $\sum_{k=0}^{{n\choose2}}w_{nk}$; see Theorem \ref{wnk}(ii) for more details.
%, which gives the total number of idempotent generating sets of $\E_n$ from $E(\D_n)$.
}
\label{tab:wnksum}
\end{table}

The next result shows that for idempotent generating sets of $\E_n$, ``minimality'' in terms of size is equivalent to ``minimality'' in terms of set containment.  

\begin{prop}\label{min_TnSn_prop}
Any idempotent generating set for $\E_n$ contains a minimal idempotent generating set.
\end{prop}

\pf The result is trivial for $n=2$, since $V=\{e_{12},e_{21}\}$ is the unique minimal (idempotent) generating set for $\E_2=\{1,e_{12},e_{21}\}$ (under either meaning of ``minimal''), so suppose $n\geq3$.  Let $U$ be an arbitrary idempotent generating set of~$\E_n$.  As mentioned above in Remark~\ref{min_rem}, $V=U\cap\D_m$ is also an idempotent generating set of $\E_n$.  It follows that the graph $\Ga=\Ga_V$ is strongly connected and complete.  Suppose $\Ga$ has ${n\choose2}+k$ edges.  If $k=0$, then $V$ is minimal, so suppose $k\geq1$, and let $i-j$ be a double edge of $\Ga$.  Let $\Ga'$ (respextively, $\Ga''$) denote the graph obtained from $\Ga$ by removing the edge $i\to j$ (respectively, $j\to i$).  It suffices to show that one of $\Ga',\Ga''$ is strongly connected.  Suppose $\Ga'$ is not strongly connected, and that  the strongly connected components of $\Ga'$ are ordered by $A_1>\cdots>A_k$.  Then, since $\Ga$ is strongly connected, it must be the case that $i\in A_k$ and $j\in A_1$.  But then, clearly, $\Ga''$ is strongly connected.  This completes the proof. \epf

\nc{\transthree}[3]{\left(\begin{smallmatrix} 1&2&3\\ #1&#2&#3 \end{smallmatrix}\right)}

\begin{rem}
A generating set for $\E_n$ need not contain a generating set of minimal size in general.  For example, the four transformations $\transthree112,\transthree131,\transthree122,\transthree232$ generate $\E_3$, although no proper subset does and $\rank(\E_3)=3$.
\end{rem}

\section{Preliminaries on $E(\TXP)$} % {\red --- is this the right name?}}

For the duration of this section, we fix integers $m,n\geq1$, the set $X=\bm\times\bn$, and the partition $\P=\{C_1,\ldots,C_m\}$ of $X$ into $m$ subsets $C_i=\{i\}\times\bn$ of size $n$.  (If $0\in\{m,n\}$, then $\TXP=\T_0$, and everything we say in this section is trivial.)  Recall that $\TXP$ is the subsemigroup of $\T_X$ consisting of all transformations $f\in \T_X$ that preserve $\P$; that is,
\[
\TXP = \set{f\in\T_X}{(\forall i\in I)(\exists j\in I) \ C_if\sub C_j}.
\]
%Note that $\TXP$ consists of only the empty mapping if $0\in\{m,n\}$.  
Note that $\TXP$ is isomorphic to $\T_n$ or $\T_m$ if $m=1$ or $n=1$, respectively.  
%Since most of the problems we investigate are either trivial in the case of $0\in\{m,n\}$, or well-known in the case of $1\in\{m,n\}$, we will typically assume that $m,n\geq2$.  However, for this section, it will be convenient to allow the possibility that $m\leq1$ or $n\leq1$.

We now describe the notation we will be using for transformations from $\TXP$.  With this in mind, let $f\in\TXP$.  There is a transformation $\fb\in\T_m$ such that, for all $i\in\bm$, $C_if\sub C_{i\fb}$.  Also, for each $i\in\bm$, there is a transformation $f_i\in\T_n$ such that $(i,j)f=(i\fb,jf_i)$ for all $j\in\bn$.  The transformation $f\in\TXP$ is uniquely determined by $f_1,\ldots,f_m\in\T_n$ and $\fb\in\T_m$, and we will write $f=[f_1,\ldots,f_m;\fb]$.  The product in $\TXP$ may easily be described in terms of this notation.  Indeed, if $f,g\in\TXP$, then $fg=[f_1g_{1\fb},\ldots,f_mg_{m\fb};\fb\gb]$.  The rule for multiplication illustrates the structure of $\TXP$ as a wreath product $\T_n\wr\T_m$, as noted in \cite{ABMS2014}.  Note that $\overline{fg}=\fb\gb$ for all $f,g\in\TXP$.

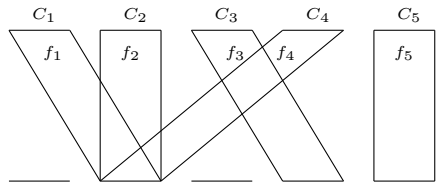
\begin{figure}[ht]
\begin{center}
\begin{tikzpicture}[xscale=.4,yscale=0.4]
\lline1
\lline2
\lline3
\lline4
\lline5
\uline1
\uline2
\uline3
\uline4
\uline5
\tll12{f_1}
\tll22{f_2}
\tll34{f_3}
\tll42{f_4}
\tll55{f_5}
%\sthicklabel11{{}_{\phantom{1}}C_1}{1.25}0
\Clab10
\Clab20
\Clab30
\Clab40
\Clab50
\end{tikzpicture}
    \caption{Diagrammatic representation of an element of $\TXP$.}
    \label{fig:pic}
   \end{center}
 \end{figure}

There is a useful way to picture a transformation $f=[f_1,\ldots,f_m;\fb]\in\TXP$.  For example, with $m=5$, and $\fb=\left(\begin{smallmatrix} 1&2&3&4&5 \\ 2&2&4&2&5 \end{smallmatrix}\right)\in\T_5$, the transformation $f=[f_1,f_2,f_3,f_4,f_5;\fb]$ is pictured in Figure \ref{fig:pic}.  

\begin{figure}[ht]
\begin{center}
\begin{tikzpicture}[xscale=.4,yscale=0.4]
\lline1
\lline2
\lline3
\lline4
\lline5
\uline1
\uline2
\uline3
\uline4
\uline5
\tll11{g_1}
\tll23{g_2}
\tll31{g_3}
\tll44{g_4}
\tll54{g_5}
\slline1{2.5}{20}
\slline2{2.5}{20}
\slline3{2.5}{20}
\slline4{2.5}{20}
\slline5{2.5}{20}
\suline1{2.5}{20}
\suline2{2.5}{20}
\suline3{2.5}{20}
\suline4{2.5}{20}
\suline5{2.5}{20}
\stll13{f_1g_2}{2.5}{20}
\stll23{f_2g_2}{2.5}{20}
\stll34{f_3g_4}{2.5}{20}
\stll43{f_4g_2}{2.5}{20}
\stll54{f_5g_5}{2.5}{20}
\slline150
\slline250
\slline350
\slline450
\slline550
\suline150
\suline250
\suline350
\suline450
\suline550
\stll12{f_1}50
\stll22{f_2}50
\stll34{f_3}50
\stll42{f_4}50
\stll55{f_5}50
\draw(20,5)node{{\large $=$}};
\Clab11
\Clab21
\Clab31
\Clab41
\Clab51
\sClab1{.5}{20}
\sClab2{.5}{20}
\sClab3{.5}{20}
\sClab4{.5}{20}
\sClab5{.5}{20}
\end{tikzpicture}
    \caption{Diagrammatic calculation of a product in $\TXP$.}
    \label{fig:prod}
   \end{center}
 \end{figure}
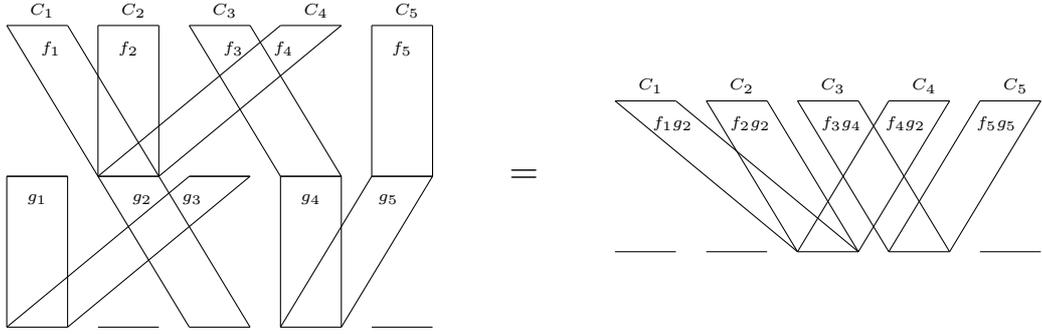

This diagrammatic representation allows for easy visualisation of the multiplication.  For example, if $f$ is as above, and if $g=[g_1,g_2,g_3,g_4,g_5;\gb]$ where $\gb=\left(\begin{smallmatrix} 1&2&3&4&5 \\ 1&3&1&4&4 \end{smallmatrix}\right)$, then the product $fg=[f_1g_2,f_2g_2,f_3g_4,f_4g_2,f_5g_5;\fb\gb]$ may be calculated as in Figure \ref{fig:prod}.

For the statement of the next result, for integers $p\geq0$ and $q\geq1$, let $\C(p,q)$ denote the set of all $q$-tuples $(a_1,\ldots,a_q)$ of non-negative integers such that $a_1+\cdots+a_q=p$.  %We write $\bl$ for the $k$-tuple $(l_1,\ldots,l_k)\in\bn^k$.

\begin{prop}\label{idempotent_prop}
A transformation $f\in\TXP$ is an idempotent if and only if
\bit
\itemit{i} $\fb\in E(\T_m)$,
\itemit{ii} $f_i\in E(\T_n)$ for all $i\in\im(\fb)$, and
\itemit{iii} $\im(f_i) \sub \im(f_{i\fb})$ for all $i\in\bm\sm\im(\fb)$.
\eit
If $m,n\geq1$, then
%The number of idempotents from $\TXP$ is equal to
\[
|E(\TXP)| = \sum_{k=1}^m \sum_{\ba} {m\choose k,a_1,\ldots,a_k} \sum_{\bl} \prod_{i=1}^k {n\choose l_i} l_i^{(a_i+1)n-l_i},
%\sum_{k=1}^m \sum_{\bl\in\bn^k} \sum_{\ba\in\C(m-k,k)} {m\choose k,a_1,\ldots,a_k} \prod_{i=1}^k {n\choose l_i} l_i^{(a_i+1)n-l_i}.
\]
where the inner sums are over all $\bl=(l_1,\ldots,l_k)\in\bn^k$ and all $\ba=(a_1,\ldots,a_k)\in\C(m-k,k)$.
\end{prop}

\pf Consider an idempotent $f\in E(\TXP)$.  Then
\[
[f_1,\ldots,f_m;\fb]=f=f^2=[f_1f_{1\fb},\ldots,f_mf_{m\fb},\fb^2],
\]
from which we immediately obtain (i).  For any $i\in\im(\fb)$, we have $i\fb=i$, so that $f_i=f_if_{i\fb}=f_i^2$, giving (ii).  Now suppose $i\in\bm\sm\im(\fb)$.  Now, $C_if\sub C_{i\fb}\cap\im(f)$.   Since $f\in E(\T_X)$, it follows that that $f$ maps $C_if$ identically onto itself.  In particular, $\{i\fb\}\times\im(f_i) = C_if\sub C_{i\fb}f = \{i\fb\}\times \im(f_{i\fb})$, establishing (iii).  Conversely, it is easy to show that any transformation $f\in\TXP$ satisfying (i--iii) is an idempotent.

In order to specify an idempotent $f=[f_1,\ldots,f_m;\fb]\in E(\TXP)$, we first specify $k=\rank(\fb)\in\bm$.  Then $\im(\fb)$  may be chosen in ${m\choose k}$ ways.  Suppose we have chosen $\im(\fb)=\{q_1,\ldots,q_k\}\sub\bm$.  For each $i\in\bk$, $f_{q_i}$ is an idempotent.  We first choose the ranks, say $l_1,\ldots,l_k\in\bn$ of these idempotents, then their images in ${n\choose l_1},\ldots,{n\choose l_k}$ ways, and then the images of the remaining elements of $C_{q_1},\ldots,C_{q_k}$ in $l_1^{n-l_1},\ldots,l_k^{n-l_k}$ ways.  We must then specify the images of the remaining elements of $X$.  To do this, we first choose the sizes of the preimages $q_1\fb^{-1},\ldots,q_k\fb^{-1}$; suppose these preimages have size $a_1+1,\ldots,a_k+1$, noting that $q_i\fb^{-1}$ contains $q_i$ for each $i$.  This may be done in ${m-k\choose a_1,\ldots,a_k}$ ways.  Once we have done this, for each $i\in\bk$, we note that each element of $\bigcup_{j\in q_i\fb^{-1}} C_j$ must map into $\im(f_{q_i})$, which has size $l_i$; so there are $l_i^{na_i}$ choices for the images of $\bigcup_{j\in q_i\fb^{-1}} C_j$.  Multiplying these terms, then adding as appropriate, and noting that ${m\choose k}{m-k\choose a_1,\ldots,a_k}={m\choose k,a_1,\ldots,a_k}$, gives the result. \epf

%\begin{rem}
%When $m=1$ in the previous result, we recover the formula for $|E(\T_n)|$ from Proposition \ref{ETn_prop}(ii).  When $n=1$, we obtain the alternative formula
%\[
%|E(\T_m)| = \sum_{k=1}^m \sum_{\ba} {m\choose k,a_1,\ldots,a_k},
%\]
%which may also be obtained directly by noting that $f\in E(\T_m)$ with $\rank(f)=k$ may be specified by first choosing $\im(f)=\{i_1,\ldots,i_k\}$ in ${m\choose k}$ ways, and then $i_1f^{-1}$ (which must include $i_1$ but none of $i_2,\ldots,i_k$) in ${m-k\choose a_1}$ ways for some $a_1\geq0$, then $i_2f^{-1}$ in ${m-k-a_1\choose a_2}$ ways for some $a_2\geq0$, and so on.  %The proposition is also (vacuously) true when $m=0$ or $n=0$.
%\end{rem}

We also give a recurrence that may be used to more easily compute the values of $|E(\TXP)|$.  For the statement of this result, it will be convenient to allow $0\in\{m,n\}$, in which case $\TXP=\T_0=E(\T_0)=\{\emptyset\}$.

\begin{prop}\label{etxp_rec}
Let $e_{mn} = |E(\TXP)|$ where $\P$ is a uniform partition of the set $X=\bm\times\bn$ into $m$ blocks of size $n$.  Then
\[
e_{0n}=1\quad\text{for all $n$},\qquad e_{mn}=\sum_{k=1}^m {m-1\choose k-1}k e_{m-k,n} \sum_{l=1}^n{n\choose l}l^{kn-l} \quad\text{for $m\geq1$.}
\]
\end{prop}

\pf It is clear that $e_{0n}=1$ for all $n$, so suppose $m\geq1$.  An idempotent $f\in E(\TXP)$ is uniquely determined by:
\ben[(i)]
\item the set $A=\set{i\in\bm}{i\fb=1\fb}$, say of size $k\in\bm$ --- there are ${m-1\choose k-1}$ choices for $A$,
\item the image $i=1\fb$, which is an element of $A$ --- there are $k$ choices for $i$, and we have $i=j\fb$ for all $j\in A$,
\item the component $f_i$, which is an idempotent from $\T_n$, say of rank $l\in\bn$ --- there are ${n\choose l}l^{n-l}$ choices for $f_i$,
\item the images of the elements of $\bigcup_{j\in A\sm\{i\}}C_j$, which must all be in $\im(f_i)$ --- there are $l^{(k-1)n}$ choices for these images, and then finally
\item the restriction of $f$ to $\bigcup_{j\in \bm\sm A}C_j$ --- there are $e_{m-k,n}$ choices for this restriction.
\een
Multiplying these values and summing over relevant $k,l$ gives the desired result. \epf

%\begin{rem}
%When $m=1$, we recover the formula for $|E(\T_n)|$ from Proposition \ref{ETn_prop}(ii).  When $n=1$, we obtain the recurrence
%\[
%|E(\T_0)|=1, \qquad |E(\T_m)| = \sum_{k=1}^m {m-1\choose k-1}  k |E(\T_{m-k})| \quad\text{for $m\geq1$}
%\]
%from \cite{DEEFHHL}.  Values of $|E(\TXP)|$ are given in Table \ref{tab:etxp}.
%\end{rem}

\begin{rem}
When $m=1$, Propositions \ref{idempotent_prop} and \ref{etxp_rec} both yield the formula for $|E(\T_n)|$ from Proposition \ref{ETn_prop}(ii).  When $n=1$, Proposition \ref{idempotent_prop} gives the alternative formula
\[
|E(\T_m)| = \sum_{k=1}^m \sum_{\ba} {m\choose k,a_1,\ldots,a_k},
\]
which may also be obtained directly by noting that $f\in E(\T_m)$ with $\rank(f)=k$ may be specified by first choosing $\im(f)=\{i_1,\ldots,i_k\}$ in ${m\choose k}$ ways, and then $i_1f^{-1}$ (which must include $i_1$ but none of $i_2,\ldots,i_k$) in ${m-k\choose a_1}$ ways for some $a_1\geq0$, then $i_2f^{-1}$ in ${m-k-a_1\choose a_2}$ ways for some $a_2\geq0$, and so on.  The $n=1$ case of Proposition \ref{etxp_rec} gives the recurrence
\[
|E(\T_0)|=1, \qquad |E(\T_m)| = \sum_{k=1}^m {m-1\choose k-1}  k |E(\T_{m-k})| \quad\text{for $m\geq1$}
\]
from \cite{DEEFHHL}.  Values of $|E(\TXP)|$ are given in Table \ref{tab:etxp}.
\end{rem}

\begin{table}[h]
\begin{center}
{
\begin{tabular}{|c||rrrrrr|}
\hline
$m\sm n$ & 0 & $1$ & $2$ & $3$ & $4$ & $5$ \\ 
\hline\hline
0   &   1   &   1   &   1   &   1   &   1   &   1   \\
1   &   1   &   1   &   3   &   10  &   41   &   196   \\
2   &   1   &   3   &   21   &   256   &   4913   &   134496   \\
3   &   1   &   10   &   189   &   9028   &   917705   &   172425016   \\
4   &   1   &   41   &   2073   &   401560   &   233777121   &   349447639616   \\
5   &   1   &   196   &   26553   &   21212980  &   74070192121   &   977698734939376   \\
\hline
\end{tabular}
}
\end{center}
\caption{Calculated values of $|E(\TXP)|$ where $\P$ is a partition of $X$ into $m$ blocks of size $n$.}
\label{tab:etxp}
\end{table}

\section{The semigroup $\EXP$}

We now move on to study the idempotent generated subsemigroup $\EXP=\la E(\TXP)\ra$ of $\TXP$.  For simplicity, we will simply write $E=E(\TXP)$ and $S=\EXP=\la E\ra$.  We will also assume from now on that $m,n\geq2$.  

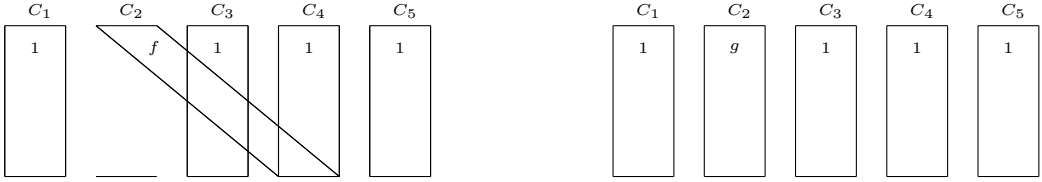
\begin{figure}[ht]
\begin{center}
\begin{tikzpicture}[xscale=.4,yscale=0.4]
\lline1
\lline2
\lline3
\lline4
\lline5
\uline1
\uline2
\uline3
\uline4
\uline5
\thickline11
\thickline24
\thickline33
\thickline44
\thickline55
\tll11{1}
\tll24{f}
\tll33{1}
\tll44{1}
\tll55{1}
\slline10{20}
\slline20{20}
\slline30{20}
\slline40{20}
\slline50{20}
\suline10{20}
\suline20{20}
\suline30{20}
\suline40{20}
\suline50{20}
\stll11{1}0{20}
\stll22{g}0{20}
\stll33{1}0{20}
\stll44{1}0{20}
\stll55{1}0{20}
\sClab10{20}
\sClab20{20}
\sClab30{20}
\sClab40{20}
\sClab50{20}
\Clab10
\Clab20
\Clab30
\Clab40
\Clab50
\end{tikzpicture}
    \caption{Diagrammatic representation of $e_{42;f}$ (left) and $g^{(2)}$ (right) from $\TXP$ with $m=5$.}
    \label{fig:e42f_f2}
   \end{center}
 \end{figure}

In what follows, certain special idempotents will play a crucial role.  With this in mind, for $i,j\in\bm$ with $i\not=j$ and for any $f\in\S_n$, we write
\[
e_{ij;f} = [1,\ldots,1,f,1,\ldots,1;e_{ij}], % \AND e_{ji;f} = [1,\ldots,1,f,1,\ldots,1;e_{ji}],
\]
where $f$ is in the $j$th position. % of $e_{ij;f}$ and in the $i$th position of $e_{ji;f}$. %{\red add picture}
Note that here $e_{ij}=\overline{e_{ij;f}}$ refers to an idempotent from $\D_m\sub\T_m$, but elsewhere we will also refer to idempotents $e_{rs}\in\D_n\sub\T_n$; however, the context should always be clear, so there should be no confusion.  Note that the transformations $e_{ij;f}$ trivially satisfy conditions (i--iii) of Proposition \ref{idempotent_prop}, so $e_{ij;f}\in E$.  %, and therefore no need to further complicate notation by writing, say $e_{ij}^{(n)}$.
If $g\in\T_n$ and $i\in\bm$, we will write $g^{(i)}=[1,\ldots,1,g,1,\ldots,1;1]$, where $g$ is in the $i$th position.
For example, with $m=5$, the transformations $e_{42;f}$ and $g^{(2)}$ are pictured in Figure \ref{fig:e42f_f2}.  For any subset $U\sub \T_n$, and for any $i\in\bm$, we write $U^{(i)} = \set{g^{(i)}}{g\in U}$.

The next result characterises the elements of $S=\la E\ra$, and also gives some information about the internal structure of $S$.  Here and elsewhere, we will write $A=B\sqcup C$ to indicate that $A$ is the disjoint union of $B$ and $C$.

\begin{prop}\label{prop:S1S2}
We have $S = S_1 \sqcup S_2$ and $S_2= S_1S_3$, where
\begin{align*}
S_1 &= \set{f\in \TXP}{\fb=1,\  f_1,\ldots,f_m\in\Sing_n},\\
S_2 &= \set{f\in\TXP}{\fb\in\T_m\sm\S_m},\\
S_3 &= \set{f\in\TXP}{\fb\in\TmSm,\ f_1,\ldots,f_m\in\S_n}.
\end{align*}
Further,
\bit
\itemit{i} $S_1=\Sing_n^{(1)}\cdots\Sing_n^{(m)}$ is the internal direct product of $m$ isomorphic copies of $\Sing_n$,
\itemit{ii} $S_2$ is isomorphic to $\T_n\wr(\TmSm)$, a wreath product of $\T_n$ with $\TmSm$, 
\itemit{iii} $S_3$ is isomorphic to $\S_n\wr(\TmSm)$, a wreath product of $\S_n$ with $\TmSm$, and
\itemit{iv} $|S| = (n^n-n!+1)^m + n^{mn}(m^m-m!)$.
\eit
% $S_1=\Sing_n^{(1)}\cdots\Sing_n^{(m)}$ is the internal direct product of $m$ isomorphic copies of $\Sing_n$, and $S_2$ is (isomorphic to) $\T_n\wr(\TmSm)$, a wreath product of $\T_n$ with $\TmSm$.  Consequently, $|S| = (n^n-n!+1)^m + n^{mn}(m^m-m!)$.
\end{prop}

\pf To prove that $S=\la E\ra\sub S_1\cup S_2$, it suffices to show that (a) $E\sub S_1\cup S_2$, and (b) $S_1\cup S_2$ is closed under multiplication on the left by the elements of $E$.  First consider an idempotent $f=[f_1,\ldots,f_m;\fb]\in E$.  Since $\fb\in E(\T_m)$, it follows that either $\fb\in\T_m\sm\S_m$ or else $\fb=1$.  In the former case, we have $f\in S_2$.  In the latter, we have %$f=f^2=[f_1^2,\ldots,f_m^2;1]$, so that 
$f_1,\ldots,f_m\in E(\T_n)\sub\Sing_n$ by Proposition~\ref{idempotent_prop}(ii), whence $f\in S_1$.  So (a) holds.  Now suppose $f\in E$ and $g\in S_1\cup S_2$.  In particular, $\fb,\gb\in\E_m$.
%We must show that $fg\in S_1\cup S_2$.  
Now, if either $\fb$ or $\gb$ belongs to $\TmSm$, then so too does $\fb\gb=\overline{fg}$, so that $fg\in S_2$.  Next, suppose $\fb=\gb=1$.  Then $f_i,g_i\in\Sing_n$ for each $i$, and it follows that $f_ig_i\in\Sing_n$ for each $i$, whence $fg=[f_1g_1,\ldots,f_mg_m;1]\in S_1$.  This completes the proof of~(b).

\begin{figure}[ht]
\begin{center}
\begin{tikzpicture}[scale=0.4]
\mfourpic13431111{0}
\mfourpic15451{h_{a_{m\text{-}2}}}11{5}
\mfourpic134111{{\raise2.5 ex\hbox{$1$}}}{\ \ {\raise0.7 ex\hbox{$1$}}}{15}
\mfourpic5345{h_{a_{1}}}111{20}
\mfourpic1344111{h_j}{25}
\mfourpic135511{h_{i}\!\!\;h_j^{\text{-}1}}1{30}
\vdotted1
\vdotted3
\vdotted4
\vdotted5
\draw(19,17.34)node{$=$};
\Clabll16{a_1}
\Clabll36{a_{m\text{-}2}}
\Clabl46{i}
\Clabl56{j}
\end{tikzpicture}
\quad
\begin{tikzpicture}[scale=0.4]
\mfourpic1344{h_{a_1}}{h_{a_{m\text{-}2}}}{h_i}{h_j}{15}
\draw[white](3,0)--(4,0); % to get previous transformation at right height
\Clabll1{3}{a_1}
\Clabll3{3}{a_{m\text{-}2}}
\Clabl4{3}{i}
\Clabl5{3}{j}
%\mfourpic124411{h_3h_4^{\!-1}}1{0}
\end{tikzpicture}
    \caption{Diagrammatic proof that $[h_1,\ldots,h_m;e_{ij}] = (e_{ji;h_ih_j^{-1}}e_{ij;h_j})\prod_{k\in\bm\sm\{i,j\}} (e_{jk;h_k}e_{kj;1})$; see the proof of Proposition \ref{prop:S1S2} for more details.}
    \label{fig:h...he}
   \end{center}
 \end{figure}
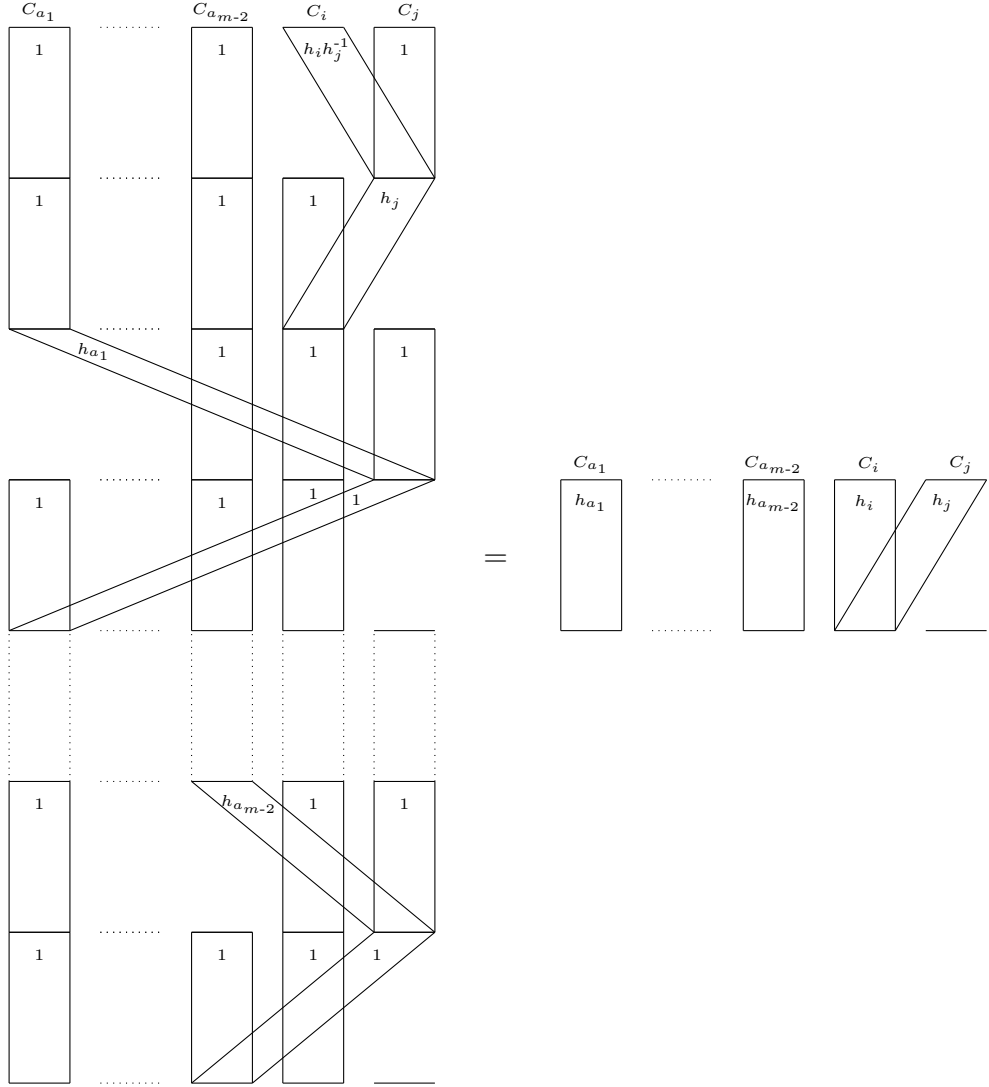

%In order to prove the reverse inclusion, we must show that $S_1$ and $S_2$ are both contained in $S=\la E\ra$. 
We now show that $S_1\cup S_2\sub S=\la E\ra$.  First, suppose $f=[f_1,\ldots,f_m;1]\in S_1$.  Then for each $i\in\bm$, $f_i\in\Sing_n=\la E(\T_n)\ra$ so that $f_i^{(i)}\in\Sing_n^{(i)}=\la E(\T_n)^{(i)}\ra\sub S$.  But then  ${f=f_1^{(1)}\cdots f_m^{(m)}\in S}$.
%First, if we pick any element $f=[f_1,\ldots,f_m;1]\in S_1$, we see that $f=f_1^{(1)}\cdots f_m^{(m)}$.  Let $U\sub E(\T_n)$ be arbitrary so that $\la U\ra=\Sing_n$.  Further, for each $i\in\bm$, we have $f_i\in\Sing_n=\la U\ra$, so that $f_i^{(i)}\in\la U^{(i)}\ra\sub\la E\ra$.  
Next suppose $f=[f_1,\ldots,f_m;\fb]\in S_2$.  First, since $\fb\in\TmSm=\la E(\D_m)\rasgp$, we write $\fb=e_{i_1j_1}e_{i_2j_2}\cdots e_{i_kj_k}$, and note that
$
f = [f_1,\ldots,f_m;e_{i_1j_1}]e_{i_2j_2;1}\cdots e_{i_kj_k;1},
%[1,\ldots,1;e_{i_2j_2}]\cdots[1,\ldots,1;e_{i_kj_k}].
$
%Since $e_{i_2j_2;1},\ldots ,e_{i_kj_k;1}\in E$, 
so it suffices to show that $[f_1,\ldots,f_m;e_{ij}]\in S$ where, for simplicity, we have written $i=i_1$ and $j=j_1$.  For each $r\in\bm$, we may write $f_r=g_rh_r$, where $g_r\in E(\T_n)$ and $h_r\in\S_n$.  Then
$
[f_1,\ldots,f_m;e_{ij}] = [g_1,\ldots,g_m;1]  [h_1,\ldots,h_m;e_{ij}].
$
Since $[g_1,\ldots,g_m;1]\in E\sub S$, it remains only to observe that
%show that $[h_1,\ldots,h_m;e_{ij}]\in S$.  But one may easily check diagramatically that
\[
[h_1,\ldots,h_m;e_{ij}] = (e_{ji;h_ih_j^{-1}}e_{ij;h_j})\prod_{k\in\bm\sm\{i,j\}} (e_{jk;h_k}e_{kj;1})\in S,
\]
where the product is calculated in ascending order of the indices $k\in\bm\sm\{i,j\}$ (although the order doesn't actually matter here).  See Figure \ref{fig:h...he} for a diagrammatic verification of this fact where, for convenience, we have written $\bm\sm\{i,j\}=\{a_1,\ldots,a_{m-2}\}$ with $a_1<\cdots<a_{m-2}$, and drawn the transformations with the blocks arranged in the order $C_{a_1},\ldots,C_{a_{m-2}},C_i,C_j$.
%an example calculation when $m=4$ and $(i,j)=(3,4)$. 

Since each $e_{ij;f}$ belongs to $S_3$, the previous paragraph also shows that $S_2=S_1S_3$.  Finally, statements (i--iv) are readily checked. \epf

Values of $|S|=|\EXP|$ are given in Table \ref{tab:exp}.

\begin{table}[ht]
\begin{center}
{
\begin{tabular}{|c||rrrrrr|}
\hline
$m\sm n$ &$0$ & $1$ & $2$ & $3$ & $4$ & $5$ \\ 
\hline\hline
0   &   1   &   1   &   1   &  1   &   1  & 1 \\
1   &   1   &   1   &   3   &  22   &   233  & 3006 \\
2   &   1   &   3   &   41   &   1942   &   185361  & 28567286\\
3   &   1   &   22   &   1371   &   423991   &   364970873     & 668031464841\\
4   &   1   &   233   &   59473   &   123528568   &   999379708193    & 22206894087218296  \\
5   &   1   &    3006    &    3077363    &    43123619167    &   3304719161323273     &     895805227489703588401   \\
\hline
\end{tabular}
}
\end{center}
\caption{Calculated values of $|\EXP|$ where $\P$ is a partition of $X$ into $m$ blocks of size~$n$.}
\label{tab:exp}
\end{table}

%the order $k=1,\ldots,i-1,i+1,\ldots,j-1,j+1,\ldots,n$.

As a consequence of the previous proof, we have the following. 

%\vspace{5mm}

\vbox{%
\begin{cor}\label{Gcor}
Let
%\[
%G_1 = \set{e_{ij}^{(k)},e_{ji}^{(k)}}{\oijn,\ k\in\bm} \quad\text{and}\quad G_2 = \set{e_{ij;f},e_{ji;f}}{\oijm,\ f\in\S_n}.
%\]
\begin{align*}
G_1 = \set{e_{ij}^{(k)},e_{ji}^{(k)}}{\oijn,\ k\in\bm} \quad\text{and}\quad G_2 = \set{e_{ij;f},e_{ji;f}}{\oijm,\ f\in\S_n}.
%G_1 &= \set{e_{ij}^{(k)},e_{ji}^{(k)}}{\oijn,\ k\in\bm},\\
%G_2 &= \set{e_{ij;f},e_{ji;f}}{\oijm,\ f\in\S_n}.
\end{align*}
Then %\nopagebreak
\begin{itemize}\begin{multicols}3
\itemit{i} $S_1=\la G_1\ra$,
\itemit{ii} $S_3=\la G_2\rasgp$, 
\itemit{iii} $S=\la G_1\cup G_2\ra$. \epf
\emc\eit
%Further, we have $S_1=\la G_1\ra$ and $S_2=\la S_1\cup G_2\ra$.
\end{cor}
}

%\begin{rem}
%It is clear that $S_1=\la G_1\ra$.  And it follows from the proof of Proposition \ref{prop:S1S2} that $$\la G_2\ra=\set{f\in\TXP}{\fb\in\TmSm,\ f_1,\ldots,f_m\in\S_n}$$ is isomorphic to a wreath product $\S_n\wr(\TmSm)$. 
%\end{rem}

The generating set $G_1\cup G_2$ from Corollary \ref{Gcor} has size $2m{n\choose 2}+2n!{m\choose2}$.  Shortly, we will see that we may reduce this set further; in fact, we will see that we can use half the elements of $G_1\cup G_2$ in the case $n\geq3$.  When $n=2$, we need all the elements of $G_1$ and half the elements of $G_2$.  %(See Theorem \ref{rank_thm}.)

\vbox{%
\begin{lemma}\label{Elemma}
Let $i\in\bm$.  Then
\bit
\itemit{i} $S\sm\Sing_n^{(i)}$ is an ideal of $S$, and  %Consequently, $S_1$ is an ideal of $S$.
\itemit{ii} any generating set for $S$ contains a generating set for $\E_n^{(i)}$.
\eit
\end{lemma}
}

\pf Part (ii) clearly follows from part (i).  To show that $S\sm\E_n^{(i)}$ is an ideal of $S$, it suffices to show that, for all $g,h\in S$, $gh\in\Sing_n^{(i)}$ implies $g,h\in\Sing_n^{(i)}$.  So suppose $g,h\in S$ are such that $gh=[g_1h_{1\gb},\ldots,g_mh_{m\gb};\gb\hb]\in\Sing_n^{(i)}$.  Since $gh\in\Sing_n^{(i)}$, it follows that $1=\gb\hb$, so that $\gb=\hb=1$.  In particular, $g,h\in S_1$ so $g_1,h_1,\ldots,g_m,h_m\in\Sing_n$.  Also, ${gh=[g_1h_1,\ldots,g_mh_m;1]\in \E_n^{(i)}}$.  So, for all $j\in\bm\sm\{i\}$, $1=g_jh_j$, giving $g_j=h_j=1$ (since also $g_j,h_j\in\E_n$).  It follows that $g,h\in\Sing_n^{(i)}$. \epf

%Recall that the \emph{relative rank} of $S$ modulo some subset $A\sub S$, denoted $\rank(S:A)$, is the least cardinality of a subset $B\sub S$ such that $S=\la A\cup B\ra$.

%\begin{lemma}
%We have $\rank(S)=\rank(S_1)+\rank(S:S_1)$ and $\rank(S:S_1)=\rank(S_3)$.
%\end{lemma}

%\pf It is clear that $S_2=S\sm S_1$ is an ideal of $S$ (because $\TmSm$ is an ideal of $\T_m$).  So it follows that $\rank(S)=\rank(S_1)+\rank(S:S_1)$.  The fact that $S=S_1\cup S_1S_3$ (and that $1\in S_1$) implies that $S=\la S_1\cup U\ra$ for any generating set $U$ of $S_3$.  This gives $\rank(S:S_1)\leq |U|\leq \rank(S_3)$.

For $\oijm$, we write $\ve_{ij}$ for the equivalence relation on $\bm$ with unique non-trivial equivalence class $\{i,j\}$.  We also write $\De = \set{(i,i)}{i\in\bm}$ for the trivial equivalence on $\bm$ (i.e., the equality relation on $\bm$).  Recall that the \emph{kernel} of a transformation $f\in\T_m$ is the equivalence $\ker(f)=\set{(i,j)\in\bm\times\bm}{if=jf}$.  Of importance is the easily checked fact that $\ker(fg)\supseteq\ker(f)$ for all $f,g\in\T_m$.

\begin{lemma}\label{Glemma}
Let $\oijm$ and $f\in\S_n$, and suppose $e_{ij;f}=gh$ where $g,h\in S$ and $g\not=1$.  Then 
\begin{itemize}\begin{multicols}{3}
\itemit{i} $\ker(\gb)=\ve_{ij}$,
\itemit{ii} $g_1,\ldots,g_m\in\S_n$, 
\itemit{iii} $g_jg_i^{-1}=f$.
\end{multicols}\eit
Consequently, if $G$ is an arbitrary generating set of $S$, then $G$ contains such an element $g$ for each $\oijm$ and $f\in\S_n$.  
%\bit
%\itemit{i} For each $i\in\bm$, $G$ contains a generating set of $\Sing_n^{(i)}$.
%\itemit{ii} 
%For each $\oijm$ and for each $f\in\S_n$, there exists some $g=[g_1,\ldots,g_m;\gb]\in G$ such that
%\begin{itemize}\begin{multicols}{3}
%\itemit{i} $\ker(\gb)=\ve_{ij}$,
%\itemit{ii} $g_1,\ldots,g_m\in\S_n$, 
%\itemit{iii} $g_jg_i^{-1}=f$.
%\end{multicols}\eit
%\eit
%Consequently, $\rank(S)\geq m\rho_n + n!{m\choose2}$, where $\rho_2=2$ and $\rho_n={n\choose2}$ if $n\geq3$.
\end{lemma}

\pf %Part (i) follows quickly from the fact, just proved, that each $S\sm\Sing_n^{(i)}$ is an ideal of $S$.  To prove part (ii), s
%Suppose $\oijm$ and $f\in\S_n$ are arbitrary.  Consider an expression $e_{ij;f} = gh_1\cdots h_k$, where $g,h_1,\ldots,h_k\in G$, and assume that $g\not=1$ (or else we delete the left-most generators $g,h_1,\ldots,h_l$ that are equal to $1$).  As usual, write $g=[g_1,\ldots,g_m;\gb]$ and, for simplicity, write $h_1\cdots h_k=p=[p_1,\ldots,p_m;\pb]$.  
Now,
$
[1,\ldots,1,f,1,\ldots,1;e_{ij}] = e_{ij;f} = gh = [g_1h_{1\gb},\ldots,g_mh_{m\gb};\gb\hb].
$
Since each $g_rh_{r\gb}$ is a permutation (either $1$ or $f$), it follows that each $g_r$ is a permutation.  Since $\E_n\cap\S_n=\{1\}$, it follows from Proposition~\ref{prop:S1S2} that either $g=[1,\ldots,1;1]=1$ or else $g\in S_2$.  We have assumed that $g\not=1$, so it follows that $g\in S_2$.  But then
$
\De \not= \ker(\gb) \sub \ker(\gb \hb) = \ker(e_{ij}) = \ve_{ij},
$
so that $\ker(\gb)=\ve_{ij}$.  It follows that $i\gb=j\gb$.  We also have $1=g_ih_{i\gb}$, so that $h_{i\gb}=g_i^{-1}$, from which it follows that
$
f=g_jh_{j\gb} = g_jh_{i\gb} = g_jg_i^{-1}.
$
So $g$ satisfies each of (i--iii).  Now suppose $G$ is an arbitrary generating set for $S$.  Then by considering an expression $e_{ij;f}=h_1\cdots h_k$, where $h_1,\ldots,h_k\in G\sm\{1\}$, we see that $h_1$ satisfies conditions (i--iii).  \epf

\begin{cor}\label{rank_cor}
We have $\rank(S)\geq m\rho_n + n!{m\choose2}$, where $\rho_2=2$ and $\rho_n={n\choose2}$ if $n\geq3$.
\end{cor}

\pf Let $G$ be an arbitrary generating set for $S$.  Lemma \ref{Elemma} (respectively, \ref{Glemma}) tells us that $G$ contains at least $m\times\rank(\E_n)=m\rho_n$ (respectively, $n!{m\choose2}$) distinct elements from $S_1$ (respectively, $S_2$).  %Lemma \ref{Elemma} tells us that $G$ contains at least $m\times\rank(\E_n)=m\rho_n$ distinct elements from $S_1$, and Lemma \ref{Glemma} tells us that $G$ contains at least another $n!{m\choose2}$ distinct elements from $S_2$.  
Since $S_1\cap S_2=\emptyset$, it follows that $|G|\geq m\rho_n+n!{m\choose2}$.  Since this is true for any generating set, it follows that this value is a lower bound for $\rank(S)$. \epf

%Next, note that there is no overlap between the elements from parts (i) and those from part (ii).  Since $\rank(\Sing_n)=\rho_n$, we see that $G$ contains at least $m\rho_n$ elements constituting the generating sets of the subsemigroups $\Sing_n^{(i)}$.  And clearly, part (ii) tells us that $G$ contains at least another $n!{m\choose2}$ elements.  This shows that any generating set of $S$ has size at least $m\rho_n+n!{m\choose2}$, so it follows that this value is a lower bound for $\rank(S)$. \epf

Our next goal is to show that the lower bound for $\rank(S)$ just obtained is in fact sharp, and also equal to $\idrank(S)$.  To do this, it is sufficient to produce any idempotent generating set of $S$ of the given size.  However, since we also aim to classify and enumerate all minimal idempotent generating sets, it will be convenient to prove a more general result.  With this in mind, we first introduce some notation.  Let
\[
\Xi = \set{(i,j)}{\oijm} \AND \Xi^{-1}=\set{(j,i)}{(i,j)\in\Xi}.
\]
If $V\sub E(\D_m)$ is such that $\TmSm = \la V\rasgp$, let
%For an idempotent generating set $V$ of $\TmSm$ with $V\sub \D_m$ (but not necessarily of minimal size), let
\[
\Xi_V = \set{(i,j)\in\Xi}{e_{ij},e_{ji}\in V}
\]
and
\[
\Phi_V = \set{(i,j)\in\Xi\cup\Xi^{-1}}{e_{ij}\in V\ \text{but}\ e_{ji}\not\in V}.
\]
Note that for any $\oijm$, exactly one of the following holds: (i) $(i,j)\in\Xi_V$, (ii) $(i,j)\in\Phi_V$, or (iii) $(j,i)\in\Phi_V$.

\begin{lemma}\label{lemS3}
Let $V\sub E(\D_m)$ with $\TmSm=\la V\rasgp$.  For each $(i,j)\in\Xi_V$, choose an ordered pair of non-empty subsets $(A_{ij},B_{ij})$ of $\S_n$ such that $\S_n=A_{ij}\sqcup B_{ij}$.  Put
\[
W = 
\bigcup_{(i,j)\in\Phi_V}\set{e_{ij;f}}{f\in\S_n}
\cup 
\bigcup_{(i,j)\in\Xi_V} \Big(
\set{e_{ij;f}}{f\in A_{ij}}\cup \set{e_{ji;f^{-1}}}{f\in B_{ij}}
\Big).
\]
Then $W\sub E$, $S_3=\la W\rasgp$ and $|W|=n!{m\choose2}$.
\end{lemma}

\begin{figure}[ht]
\begin{center}
\begin{tikzpicture}[scale=0.35]
\mtwopicl111102
\mtwopic22g101
\draw(10,9.84)node{$=$};
\mtwopicl22gg9{1.5}
\mtwopicl22gg{21}3
\mtwopic22gg{21}2
\mtwopic22gg{21}0
\draw[dotted](24,5)--(24,10);
\draw[dotted](26,5)--(26,10);
\draw[dotted](27,5)--(27,10);
\draw[dotted](29,5)--(29,10);
\draw(31,9.84)node{$=$};
\mtwopicl22{g^q}{g^q}{30}{1.5}
\draw(40,9.84)node{$=$};
\mtwopicl2211{39}{1.5}
%\draw(19,14.84)node{$=$};
\end{tikzpicture}
    \caption{Diagrammatic proof that $(e_{ij;1}e_{ji;g})^q=e_{ji;1}$;
    %$(e_{ij;f^{-1}}e_{ji;1})^q=e_{ji;f}$; 
    see the proof of Lemma \ref{lemS3} for more details.}
    \label{fig:gq}
   \end{center}
 \end{figure}
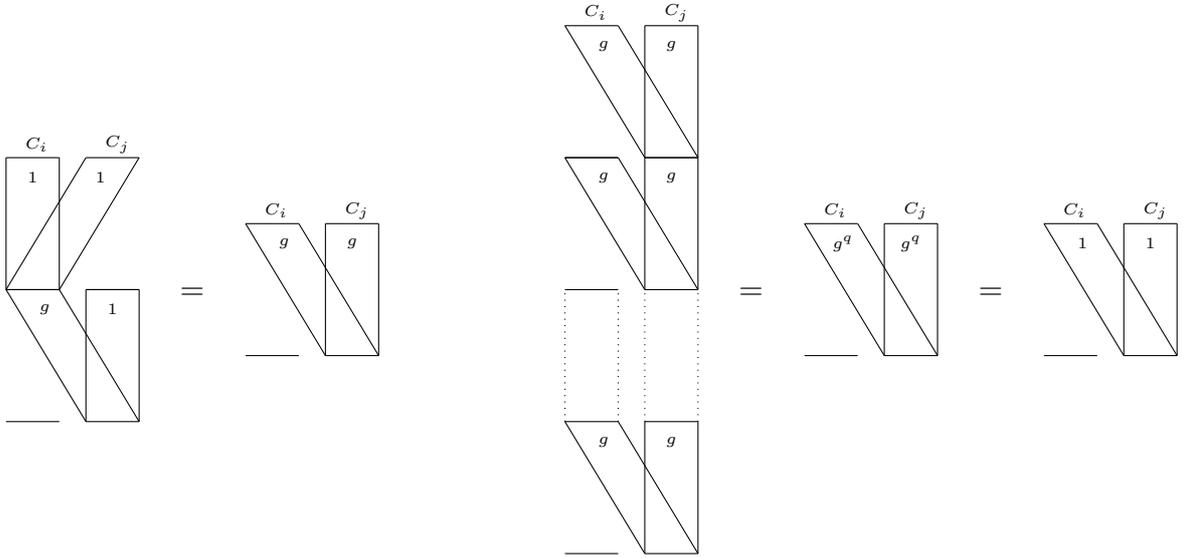

\pf It is clear that $W\sub E$ has the required size. So, by Corollary \ref{Gcor}(ii), it remains to show that $\la W\rasgp$ contains each $e_{ij;f}\in G_2$.  To simplify notation, let $T=\la W\rasgp$.  
We first claim that $T$ contains $e_{ij;1}$ for each $e_{ij}\in V$.  Now $e_{ij;1}$ belongs to $W\sub T$ for each $(i,j)\in \Phi_V$.  Next, suppose $(i,j)\in\Xi_V$.  Then either (i) $e_{ij;1}\in W$ or (ii) $e_{ji;1}\in W$.  Suppose~(i) holds.  We will show that $e_{ji;1}\in T$.  (A similar argument shows that (ii) implies $e_{ij;1}\in T$.)  Since $B_{ij}\not=\emptyset$, $W$ contains some $e_{ji;g}$ with $g\in \S_n$.  Let $q\geq1$ be such that $g^q=1$.  Then $e_{ji;1}=(e_{ij;1}e_{ji;g})^q\in T$, as we show in Figure~\ref{fig:gq} where, for convenience, we have only pictured the action of the transformations on the blocks $C_i$ and $C_j$, all other blocks being mapped identically.
This completes the proof of the claim.  Since $\TmSm=\la V\rasgp$, every $e_{ij}$ is a product of elements from $V$, so it also follows that~$T$ contains every $e_{ij;1}$.

Next, we show that $T$ contains each $e_{ji;f}$ with $(i,j)\in\Phi_V$.  Now $W$ contains $e_{ij;f^{-1}}$ and we know that $e_{ji;1}\in T$.  %We show in Figure \ref{fig:fq} 
One may check diagrammatically that $e_{ji;f}=(e_{ij;f^{-1}}e_{ji;1})^q\in T$, where $q\geq1$ is such that $f^q=1$.
%
%So it follows that $T$ contains $e_{ij;f^{-1}}e_{ji;1}=[1,f^{-1};e_{ji}]$.  (Here, again, we are only listing the $i$th and $j$th components, but we are writing the $i$th component first, even if $j<i$.)  Let $q\geq1$ be such that $f^q=1$.  Then $T$ contains $[1,f^{-1};e_{ji}]^q=[f^{1-q},f^{-q};e_{ji}]=[f,1;e_{ji}]=e_{ji;f}$.
%
Now let $(i,j)\in\Xi_V$.  Then $W$ contains each $e_{ij;f}$ with $f\in A_{ij}$, and each $e_{ji;f^{-1}}$ with $f\in B_{ij}$.  The proof of the lemma will be complete if we can show that $T$ contains (i) each $e_{ij;f}$ with $f\in B_{ij}$, and (ii) each $e_{ji;f^{-1}}$ with $f\in A_{ij}$.  We will show that (i) is true, with (ii) being similar.  Let $f\in B_{ij}$.  Then $e_{ji;f^{-1}}\in W$, and we also know that $e_{ij;1}\in T$.  As above, we see that $e_{ij;f} = (e_{ji;f^{-1}}e_{ij;1})^q$ for some $q\geq1$. \epf

%\begin{figure}[ht]
%\begin{center}
%\begin{tikzpicture}[scale=0.35]
%\mtwopicl111{f^{-1}}02
%\mtwopic221101
%\draw(10,9.84)node{$=$};
%\mtwopicl221{f^{-1}}9{1.5}
%%
%\mtwopicl221{f^{-1}}{21}3
%\mtwopic221{f^{-1}}{21}2
%\mtwopic221{f^{-1}}{21}0
%\draw[dotted](24,5)--(24,10);
%\draw[dotted](26,5)--(26,10);
%\draw[dotted](27,5)--(27,10);
%\draw[dotted](29,5)--(29,10);
%\draw(31,9.84)node{$=$};
%\mtwopicl22{f^{1-q}}{f^{-q}}{30}{1.5}
%\draw(40,9.84)node{$=$};
%\mtwopicl22f1{39}{1.5}
%%\draw(19,14.84)node{$=$};
%\end{tikzpicture}
%    \caption{Diagrammatic proof that %$(e_{ij;1}e_{ji;g})^q=e_{ji;1}$
%    $(e_{ij;f^{-1}}e_{ji;1})^q=e_{ji;f}$; 
%    see the proof of Lemma \ref{lemS3} for more details.}
%    \label{fig:fq}
%   \end{center}
% \end{figure}

We are now ready to prove one of the main results of the paper.

\begin{thm}\label{rank_thm}
Let $U_1,\ldots,U_m$ be minimal idempotent generating sets of $\E_n$, and let $W$ be as in Lemma \ref{lemS3}.  Then
$
U_1^{(1)}\cup\cdots\cup U_m^{(m)} \cup W
$
is an idempotent generating set of $S$ of size $m\rho_n + n!{m\choose2}$.  Consequently, $$\rank(S)=\idrank(S) = m\rho_n + n!{m\choose2},$$ where $\rho_2=2$ and $\rho_n={n\choose2}$ if $n\geq3$.  
\end{thm}

\pf Put $V=U_1^{(1)}\cup\cdots\cup U_m^{(m)}\cup W$.  Since $V\sub E$ and $|V|=m\rho_n + n!{m\choose2}$, and since we already know that $\idrank(S)\geq\rank(S)\geq m\rho_n+n!{m\choose2}$, by Corollary \ref{rank_cor}, it suffices to show that $\la V\ra=S$.  But $S=S_1\cup S_1S_3$ by Proposition \ref{prop:S1S2}, and clearly $S_1=\la U_1^{(1)}\cup\cdots\cup U_m^{(m)} \ra$.  Since $\la W\rasgp=S_3$ by Lemma \ref{lemS3} (and since $1\in S_1$), the proof is complete. \epf

\begin{rem}
With one exception, the expression for $\rank(S)$ given in Theorem \ref{rank_thm} is valid if $1\in\{m,n\}$, if we also define $\rho_1=\rank(\E_1)=\idrank(\E_1)=0$.  When $m=1$, the formula reduces to $\rho_n$.  When $n=1$, it reduces to ${m\choose2}$.  So in both cases, we see that it agrees with Theorem~\ref{Howie1978_thm}, unless $(m,n)=(2,1)$, where the formula gives ${2\choose2}=1$, even though ${\rank(S)=\rank(\E_2)=2}$.
The expression is also valid if $m=0$, giving a value of $0$.  But it is not valid for $n=0$ unless $m\leq1$ (even if we define $\rho_0=0$), as it gives ${m\choose2}$ rather than $0$.
%
%If we further define $\rho_n=\rank(\E_n)={n\choose2}=0$ for $n\leq1$, then the formula for $\rank(S)=\idrank(S)$ given in Theorem \ref{rank_thm} is valid for $m\leq1$ or $n\leq1$ except for the case in which $(m,n)=(2,1)$.
%
Calculated values of $\idrank(S)=\rank(S)$ are given in Table \ref{tab:txp}.
\end{rem}

%As a corollary of Lemmas \ref{Glemma} and \ref{lem1}, we immediately obtain the following.

%\begin{thm}\label{rank_thm}
%We have $\rank(S)=\idrank(S) = m\rho_n + n!{m\choose2}$, where $\rho_2=2$ and $\rho_n={n\choose2}$ if $n\geq3$. \epfres
%\end{thm}

\begin{table}[h]
\begin{center}
{
\begin{tabular}{|c||rrrrrrrrrr|}
\hline
$m\sm n$ & 1 & 2 & 3 & 4 & 5 & 6 & 7 & 8 & 9 & 10 \\ 
\hline\hline
1     &      0     &      2    &       3     &      6       &   10         & 15           &    21           &28                &36               &45\\
2     &      2     &      6    &      12    &      36     &    140      &   750       &   5082        &40376          &362952       &3628890\\
3     &      3     &     12   &       27   &       90    &     390     &   2205     &  15183       &121044        &1088748     & 10886535\\
4     &      6     &     20   &       48   &      168   &      760    &    4380    &   30324      &242032        &2177424     & 21772980\\
5     &     10    &      30  &        75  &       270  &      1250  &      7275  &     50505    &  403340      & 3628980    &  36288225\\
6     &     15    &      42  &       108 &        396 &       1860 &      10890&       75726  &    604968    & 5443416    &  54432270\\
7     &     21    &      56  &       147 &        546 &       2590 &      15225&      105987 &     846916   &  7620732   & 76205115\\
8     &     28    &      72  &       192 &        720 &       3440 &      20280&      141288 &    1129184  &  10160928 &  101606760\\
9     &     36    &      90  &       243 &        918 &       4410 &      26055&      181629 &    1451772  &  13064004 &  130637205\\
10 &         45  &     110 &        300&        1140&        5500&       32550&      227010&     1814680 &   16329960&   163296450\\
\hline
\end{tabular}
}
\end{center}
\caption{Calculated values of $\rank(\EXP)=\idrank(\EXP)$ where $\P$ is a partition of $X$ into $m$ blocks of size~$n$.}
\label{tab:txp}
\end{table}

Next we aim to enumerate the minimal idempotent generating sets of $S$.  Part of doing this involves showing that every such generating set has the form described in the previous theorem.
%Actually, we have already done half the work in Lemma \ref{lem1}.
We also aim to show that any idempotent generating set of $S$ contains a minimal (idempotent) generating set.  In order to accomplish these aims simultaneously, we will need to prove a number of technical lemmas.

\begin{lemma}\label{newlemma1}
Let $U\sub E$ and $1\leq r<s\leq m$.  Suppose $e_{rs;1}=g_1\cdots g_k$, where $g_1,\ldots,g_k\in U$ and $k$ is minimal among all such expressions of $e_{rs;1}$ as a product of elements from $U$.  Then $g_1,\ldots,g_k\in G_2=\set{e_{ij;f},e_{ji;f}}{\oijm,\ f\in\S_n}$.
\end{lemma}

\pf Write $g_i=[g_{i1},\ldots,g_{im};\gb_i]$ for each $i$.  
%Since $e_{rs}=\gb_1\cdots\gb_k$, it follows that $\rank(\gb_i)\geq n-1$ for all $i$.  
First note that $\ve_{rs}=\ker(e_{rs})=\ker(\gb_1\cdots\gb_k)\supseteq\ker(\gb_1)$, so that $\ker(\gb_1)=\ve_{rs}$ or $\De$, the trivial equivalence.  As in the proof of Lemma \ref{Glemma}, $\ker(\gb_1)=\De$ would imply that $g_1=1$, so that $e_{rs;1}=g_2\cdots g_k$, contradicting the minimality of~$k$.  So we must in fact have $\ker(\gb_1)=\ve_{rs}$.  Since $\gb_1$ is an idempotent, it follows that $\gb_1=e_{rs}$ or $e_{sr}$.  Suppose $\gb_1=e_{rs}$, the other case being similar.  Since $g_1\cdots g_k=e_{rs;1}$ is injective when restricted to each $C_i$, we see that $g_{1i}\in\S_n$ for all $i\in\bm$.  Since $g_1$ is an idempotent, Proposition~\ref{idempotent_prop} gives $g_{1i}\in E(\T_n)$ for all $i\in\im(\gb_1)=\bm\sm\{s\}$, so that $g_{1i}=1$ for all such $i$.  It follows that $g_1=e_{rs;g_{1s}}\in G_2$.  (Further considerations reveal that in fact $g_{1s}=1$, but we do not need to know this.)

As an inductive hypothesis, suppose $g_1,\ldots,g_{l-1}\in G_2$ for some $2\leq l\leq k$.  We must show that $g_l\in G_2$.
%Now suppose $g_l\not\in G_2$ for some $2\leq l\leq k$, and suppose $l$ is the least such index.  
Put $h=g_1\cdots g_{l-1}$ and write $h=[h_1,\ldots,h_m;\hb]$.  Again, we have $h_1,\ldots,h_m\in\S_n$.  Note that $n-1=\rank(e_{rs})=\rank(\hb\gb_l\cdots\gb_k)\leq\rank(\hb)\leq\rank(\gb_1)=n-1$, so it follows that $\rank(\hb)=n-1$.  Say $\im(\hb)=\bm\sm\{q\}$.  We also have $\rank(\gb_l)\geq n-1$.  Suppose first that $\rank(\gb_l)=n$, so that $\gb_l=1$.  Then Proposition \ref{idempotent_prop} gives $g_{li}\in E(\T_n)$ for each $i$.  If $g_{li}\not=1$ for some $i\in\bm\sm\{q\}$, then for any $j\in\bm$ with $j\hb=i$, the restriction of $hg_l$ to $C_j$ would not be injective, and hence neither would the restriction of $e_{rs;1}$ to $C_j$, a contradiction.  So it follows that $g_{li}=1$ for all $i\in\bm\sm\{q\}$.  But then $g_l$ acts as the identity on $\bigcup_{i\in\bm\sm\{q\}}C_i=\im(h)$, and so $hg_l=h$, giving $e_{rs;1}=g_1\cdots g_{l-1}g_{l+1}\cdots g_k$, again contradicting the minimality of $k$.  Now suppose $\rank(\gb_l)=n-1$.  Since $\rank(\hb\gb_l)=n-1$, and since $\im(\hb)=\bm\sm\{q\}$, it follows that $\gb_l=e_{pq}$ or $e_{qp}$ for some $p\in\bm\sm\{q\}$.  A similar argument to that just used shows that if $\gb_l=e_{pq}$, then $hg_l=h$, and we again obtain a contradiction.
%Suppose first that $\gb_l=e_{pq}$.  Then $g_{li}\in E(\T_n)$ for all $i\in\im(e_{pq})=\bm\sm\{q\}$.  But also $g_l$ must be injective when restricted to $\im(h)=\bigcup_{i\in\bm\sm\{q\}}C_i$, so it follows that $g_{li}=1$ for all $i\in\bm\sm\{q\}$, so that $hg_l=h$, again leading to a contradiction.  
Finally, suppose $\gb_l=e_{qp}$.  We again have $g_{li}\in E(\T_n)$ for all $i\in\bm\sm\{p\}$, and $g_{lj}\in\S_n$ for all $j\in\bm\sm\{q\}$.  
%Since $g_l$ must be injective on $\im(h)$, 
It follows that $g_{li}=1$ for $i\in\bm\sm\{p,q\}$.  
We also have %$g_{lq}\in E(\T_n)$, $g_{lp}\in\S_n$, and 
$\im(g_{lq})\supseteq \im(g_{lp})=\bn$, by Proposition \ref{idempotent_prop}, so $g_{lq}=1$.  It follows that $g_l=e_{qp;g_{lp}}\in G_2$, completing the proof. \epf

\begin{lemma}\label{newlemma2}
Suppose $W\sub E$ and $S=\la W\ra$.  Then there exist subsets $W_1,W_2\sub W$ such that
\bit
%\begin{itemize}\begin{multicols}3
\itemit{i} $|W_1|=m\rho_n$ and  $S_1=\la W_1\ra$, 
\itemit{ii} $|W_2|=n!{m\choose2}$ and
\begin{itemize}
%\end{multicols}
%\itemit{iv} $W_2\sub G_2=\set{e_{ij;f},e_{ji;f}}{\oijm,\ f\in\S_n}$,
\itemit{a} for all $\oijm$ and for all $f\in\S_n$, either $e_{ij;f}\in W_2$ or $e_{ji;f^{-1}}\in W_2$,
\itemit{b} the set $V=\set{e_{ij}}{(\exists f\in\S_n)\ e_{ij;f}\in W_2}$ generates $\TmSm$ (as a semigroup).
\eit
\eit
%Further, we have $S_1=\la W_1\ra$, $S_3=\la W_2\ra$, and $S=\la W_1\cup W_2\ra$.
%\bit
%%\begin{itemize}\begin{multicols}3
%\itemit{i} $|W_1|=m\rho_n$,
%\itemit{ii} $S_1=\la W_1\ra$, 
%\itemit{iii} $|W_2|=n!{m\choose2}$,
%%\end{multicols}
%%\itemit{iv} $W_2\sub G_2=\set{e_{ij;f},e_{ji;f}}{\oijm,\ f\in\S_n}$,
%\itemit{iv} for all $\oijm$ and for all $f\in\S_n$, either $e_{ij;f}\in W_2$ or $e_{ji;f}\in W_2$,
%\itemit{v} the set $V=\set{e_{ij}}{(\exists f\in\S_n)\ e_{ij;f}\in W_2}$ generates $\TmSm$.
%\eit
\end{lemma}

\pf By Lemma~\ref{Elemma}(ii), $W$ contains a generating set $U_i$ for each $\Sing_n^{(i)}$, and we clearly have $U_i\sub E$ for each $i$.  By Proposition \ref{min_TnSn_prop}, each $U_i$ contains a minimal idempotent generating set, $U_i'$, of $\Sing_n^{(i)}$.  Put $W_1=U_1'\cup\cdots\cup U_m'$, noting that $|W_1|=m\rho_n$ and $\la W_1\ra=S_1$, establishing (i).

Now we turn to (ii).  First note that by Lemma \ref{Glemma}, $W$ contains a subset $W_2'$ with $|W_2'|=n!{m\choose2}$ and such that, for all $\oijm$ and $f\in\S_n$, there exists some $g=[g_1,\ldots,g_m,\gb]\in W_2'$ with
\bit\bmc{3}
\item[(1)]  $\ker(\gb)=\ve_{ij}$,
\item[(2)] $g_1,\ldots,g_m\in\S_n$,
\item[(3)] $g_jg_i^{-1}=f$.
\end{multicols}\eit
Consider now some $g\in W_2'$ satisfying (1--3) above.  Since $W$ consists entirely of idempotents, we immediately have $\gb=e_{ij}$ or $\gb=e_{ji}$.  In the former case, since $g$ is an idempotent, we see that $g_k=1$ for all $k\in\bm\sm\{j\}$, and also $f=g_jg_i^{-1}=g_j$, so that $g=e_{ij;f}$.  Similarly, in the latter case, we deduce that $g=e_{ji;f^{-1}}$.  It follows that $W_2'$ satisfies (a).  Now put $$V'=\set{e_{ij}}{(\exists f\in\S_n)\ e_{ij;f}\in W_2'},$$ and consider the digraph $\Ga=\Ga_{V'}$, as defined in Theorem~\ref{Howie1978_thm}.  By (a), and the definition of $V'$, we see that $\Ga$ is complete.  Let
\[
\Psi=\Psi_{V'}=\set{(i,j)\in\bm\times\bm}{\text{$i\not=j$ and there is no path in $\Ga$ from $i$ to $j$}}.
\]
%and put $\psi=\psi_{V'}=|\Psi|$.  
If $|\Psi|=0$, then $\Ga$ is strongly connected, from which it follows from Theorem~\ref{Howie1978_thm} that $V'$ generates $\TmSm$, showing that $W_2'$ satisfies (b), and completing the proof.  So suppose instead that $|\Psi|\geq1$, and let $(r,s)\in\Psi$.  Consider an expression $e_{rs;1} = g_1\cdots g_k$, where $g_1,\ldots,g_k\in W$ and $k$ is minimal.  By Lemma \ref{newlemma1}, each of the $g_i$ belongs to $G_2$, say $g_i=e_{a_ib_i;f_i}$.  It then follows that $e_{rs}=\gb_1\cdots\gb_k=e_{a_1b_1}\cdots e_{a_kb_k}$.  But since $(r,s)\in\Psi$, it follows that not all of $e_{a_1b_1},\ldots,e_{a_kb_k}$ belong to $V'$.  Suppose $p_1,\ldots,p_l\in\bk$ are such that $e_{a_{p_1}b_{p_1}},\ldots,e_{a_{p_l}b_{p_l}}\not\in V'$, but the other $e_{a_ib_i}$ do belong to $V'$.  Now, $W_2'$ contains each $e_{b_{p_i}a_{p_i};f_{p_i}^{-1}}$.  Put
\[
W_2'' = \big(W_2'\sm\{e_{b_{p_1}a_{p_1};f_{p_1}^{-1}},\ldots,e_{b_{p_l}a_{p_l};f_{p_l}^{-1}}\}\big)\cup\{e_{a_{p_1}b_{p_1};f_{p_1}},\ldots,e_{a_{p_l}b_{p_l};f_{p_l}}\}.
\]
Then $W_2''\sub W$ certainly still satisfies condition (a).  Analogously to $V'$ and $\Psi_{V'}$, we define $V''$ and $\Psi_{V''}$, and see that
%And the set
%\[
%\Psi_{V''}=\set{(i,j)\in\bm\times\bm}{\text{$i\not=j$ and there is no path in $\Ga_{V''}$ from $i$ to $j$}}
%\]
%satisfies 
$|\Psi_{V''}|<|\Psi_{V'}|$.  Continuing in this fashion, we will eventually arive at a subset $W_2\sub W$ satisfying both (a) and (b). \epf %Note that $W_1\cup W_2$ has the form described in Lemma \ref{lemS3}, so it follows that $S=\la W_1\cup W_2\ra$.  This completes the proof. \epf

\newpage

\begin{thm}\label{enum_thm}
\bit
\itemit{i} Every idempotent generating set of $S$ contains a minimal idempotent generating set.
\itemit{ii} Every minimal idempotent generating set of $S$ is of the form described in Theorem \ref{rank_thm}.
\itemit{iii} The number of minimal idempotent generating sets of $S$ is equal to
\[
\si_n^m\times \sum_{k=0}^{{m\choose 2}} w_{mk}(2^{n!}-2)^k,
\]
where $\si_2=1$ and $\si_n=w_n$ otherwise.  Recurrences for the numbers $w_n$ and $w_{nk}$ are given in Theorems \ref{wn} and \ref{wnk}.
\eit
\end{thm}

\pf Part (i) follows from Lemma \ref{newlemma2} and Theorem \ref{rank_thm}.  For part (ii), suppose $W$ is a minimal idempotent generating set of $S$.  Then we must have $W=W_1\cup W_2$ where $W_1,W_2$ are as described in Lemma~\ref{newlemma2}, so it follows that $W$ has the form described in Theorem \ref{rank_thm}.  By part (ii), a minimal idempotent generating set for $S$ is completely determined by (using the notation of Lemma \ref{lemS3} and Theorem \ref{rank_thm}):

~\qquad (a) $V$, \qquad (b) $U_1,\ldots,U_m$, \qquad (c) the subsets $(A_{ij},B_{ij})$ for each $(i,j)\in\Xi_V$.

%\bit\begin{multicols}{2}
%\item[(i)] $V$,
%\item[(ii)] $U_1,\ldots,U_m$,
%\item[(iii)] the subsets $(A_{ij},B_{ij})$ for each $(i,j)\in\Xi_V$.
%\end{multicols}\eit
There are $\si_n^m$ ways to choose the subsets $U_1,\ldots,U_m$.  For each $0\leq k\leq{m\choose2}$, there are $w_{mk}$ choices of $V$ with $|\Xi_V|=k$ and, for such a $V$, there are $2^{n!}-2$ ways to choose the subsets $(A_{ij},B_{ij})$ for each $(i,j)\in\Xi_V$. \epf

\begin{rem}
With one exception, the expression given in Theorem \ref{enum_thm} is valid if $1\in\{m,n\}$, if we also define $\si_1=w_1=1$.  When $m=1$, the formula reduces to $\si_n$ (assuming $0^0=1$).  When $n=1$, it reduces to $w_m$.  So in both cases, we see that it agrees with Theorem \ref{wn}, unless $(m,n)=(2,1)$, where the formula gives $w_2=0$, even though there is actually a unique (minimal idempotent) generating set of $\EXP\cong\E_2$.
Some calculated values are given in Table \ref{tab:minexp}.
%Some calculated values of $\si_n^m\times \sum_{k=0}^{{m\choose 2}} w_{mk}(2^{n!}-2)^k$ are given in Table \ref{tab:minexp}.
\end{rem}

\begin{table}[h]
\begin{center}
{
\begin{tabular}{|c||rrrr|}
\hline
$m\sm n$ & 1 & 2 & 3 & 4  \\ 
\hline\hline
1     &   1   &  1    &  2    &   24   \\
2     &  1    &  2    &  248    &    9663675264  \\
3     &  2    &  46    &   2094128   &  65281994259188583864812544    \\
4     &   24   &   3608   &   1099477716608   &    $7.398852038987696\times 10^{48}$  \\
\hline
\end{tabular}
}
\end{center}
\caption{The number of minimal idempotent generating sets of $\EXP$ where $\P$ is a partition of $X$ into $m$ blocks of size~$n$.}
 %\caption{Calculated values of $\si_n^m\times \sum_{k=0}^{{m\choose 2}} w_{mk}(2^{n!}-2)^k$; see Theorem \ref{enum_thm} for more details.}
\label{tab:minexp}
\end{table}

\section*{Acknowledgement}

The authors wish to thank the referee for their careful reading and helpful comments that led to increased clarity and truth.

\footnotesize
\def\bibspacing{-1.1pt}
\bibliography{biblio}
\bibliographystyle{plain}
\end{document}